\documentclass[14]{article}
\usepackage[margin=1in]{geometry}
\usepackage{amsfonts,amsmath,amssymb}
\numberwithin{equation}{section}

\title{On Musielak $N$-functions}
\author{Abdulhameed Qahtan Abbood Altai \\ \small
Department of Mathematics and Statistics, University of Maryland, Baltimore County\\
\small 1000 Hilltop Circle\\  \small Baltimore, MD 21250, \\ \small altai1@umbc.edu}

\begin{document}

\maketitle

\paragraph{Abstract.} 
In this paper, the concept of Musielak $N$-functions and Musielak-Orlicz spaces generated by them well be introduced. Facts and results of the measure theory will be applied to consider properties, calculus and basic approximation of Musielak $N$-functions and their Musielak-Orlicz spaces. Finally, the relationship between Musielak $N$-functions and Musielak-Orlicz functions and thier Musielak-Orlicz spaces will be considered using facts and results of the measure theory too.


\section{Introduction} 
$N$-functions, Orlicz functions and Orlicz classes and Orlicz spaces generated by $N$-functions and Orlicz functions have been studied by many mathematicians as in [24],[17],[28],[18],[25],[3],[4],[5],[26],[11],[20],[9]. Musielak-Orlicz functions and Musielak-Orlicz spaces generated by Musielak-Orlicz functions have been originated and developed by [23],[22],[21] where $f \in L_{MO}\left( \Omega,\Sigma, \mu \right)$ if and only if $ \int_{\Omega} MO(t,f(t)) d\mu < \infty$. Their properties have been studied by [13],[16],[8],[14],[15],[16],[33] and their applications can be found in differential equations [7],[10], fluid dynamics [29],[31], statistical physics[1],integral equations [17], image processing [2],[6],[12] and many other applications [27]. So, because such increasingly importance to these concepts in the modeling of modern materials, we want to investigate properties, calculus and basic approximations of Musielak $N$-functions and  Musielak-Orlicz functions and their Musielak-Orlicz spaces using the measure theory where this will help us to consider $\mu-$almost everywhere property, supremum, infimum, limit, convergence and basic convergence of Musielak $N$-functions, Musielak-Orlicz functions and Musielak-Orlicz spaces generated by them by functioning facts and results of the measure theory and getting advantages from that to consider these concepts. The concept of Musielak $N$-function $M(t,u)$ is a generalization to the concept of $N$-function $M(u)$, where $M(t,u)$ may vary with location in space, whereas the Musielak-Orlicz function $MO(t,u)$ is a generalization to the concept of Orlicz functions $O(u)$, where $MO(t,u)$ may vary with location in space. The Musielak-Orlicz function $MO(t,u)$ is defined on $\Omega \times [0,\infty)$ into $[0,\infty)$ where for $\mu-$a.e. $t \in \Omega, MO(t,.)$ is Orlicz function of $u$ on $[0,\infty)$ and for each $u \in [0,\infty), MO(.,u)$ is a $\mu-$measurable function of $t$ on $\Omega$ and $(\Omega,\Sigma,\mu)$ is a measure space [3],[15]. So, we are going to define the Musielak $N$-function $M(t,u)$ on $\Omega \times \mathbb{R}$ into $\mathbb{R}$ in similar way, that is, for $\mu-$a.e. $t \in \Omega, M(t,.)$ is $N$-function of $u$ on $\mathbb{R}$ and for each $u \in \mathbb{R}, M(.,u)$ is a $\mu-$measurable function of $t$ on $\Omega$ and $(\Omega,\Sigma,\mu)$ is a measure space. The novelty to define the Musielak $N$-function $M(t,u)$ by this way is to get benefits from the results of the measure theory and use them to consider properties, calculus, and basic convergence of Musielak $N$-functions and Musielak-Orlicz spaces generated by them and the relationship between Musielak $N$-functions and Musielak-Orlicz functions and their Musielak-Orlicz spaces generated by them where this will give us more flexibility to pick a suitable measurable set $\Omega$ and then the functional $ \int_{\Omega}M(t,\Vert f(t) \Vert_{BS})d\mu $ defined on it as we will see in section 3. So, the paper is organized as follows. Definition of Musielak $N$-function, developing preliminaries results about the equivalent definition of Musielak $N$-function and studying continuity of Musielak $N$-function are introduced in section 2. Definition of Musielak-Orlicz space generated by a Musielak $N$-function, and using facts and results of the measure theory to study properties, calculus and basic approximation of Musielak $N$-functions and Musielak-Orlicz spaces generated by them in section 3. The relationship between Musielak $N$-functions and Musielak-Orlicz functions and Musielak-Orlicz spaces generated by them respectively using facts and results of the measure theory are introduced also in section 4. Examples of Musielak $N$-functions and Musielak-Orlicz functions that are not Musielak $N$-functions will be in section 5. The conclusion will be in section 6.

\section{Preliminary Results}
In this section, we introduce the concept of the Musielak $N$-function and some results about the equivalent definition and the continuity of the Musielak $N$-functions.

\paragraph{Defination 2.1 (Musielak $N$-function).}

Let $(\Omega, \Sigma, \mu)$ be a measure space. A function $M:\Omega \times \mathbb{R} \to \mathbb{R}$ is called a Musielak $N$-function if:
\begin{enumerate}
\item for $\mu-$a.e. $t\in \Omega, M(t,u)$ is even convex of $u$ on $\mathbb{R}$
\item for $\mu-$a.e. $t\in \Omega, M(t,u)> 0$ for any $u>0$
\item for $\mu-$a.e. $t\in \Omega, \lim \limits_{u \to 0} \frac{M(t,u)}{u} = 0$
\item for $\mu-$a.e. $t\in \Omega, \lim \limits_{u \to + \infty} \frac{M(t,u)}{u} = + \infty$
\item for each $u \in \mathbb{R}, M(t,u)$ is a $\mu-$measurable function of $t$ on $\Omega$.
\end{enumerate}

The conditions $(1-4)$ guarantee that for $\mu-$a.e. $t\in \Omega, M(t,.)$ is $N$-function of $u$ on $\mathbb{R}$ and the condition $(5)$ guarantees the measurability of $M(.,u)$ of $t$ on $\Omega.$ \\

The following theorem is a generalization to theorem 1.1[17] which is necessary to proof the next theorem where it gives us the equivalent definition of Musielak $N$-function.
 
\paragraph{Theorem 2.1.}For $\mu-$a.e. $t\in \Omega$, every convex function $M(t,u):\Omega \times \left[ a,b \right] \longrightarrow \mathbb{R}$ of $u$ on $\left[ a,b\right]$ which satisfies the condition $M(t,a)=0$ can be represented in the form 
\begin{align}
M(t,u) = \int_{a}^{u} p(t,s)ds,
\end{align}
where $p(t,u):\Omega \times \left[ a,b \right] \longrightarrow \mathbb{R}$ is, for $\mu-$a.e. $t\in \Omega$, a non-decreasing, right-continuous function of $u$ on $\left[ a,b \right], (\Omega, \Sigma, \mu)$ is a measure space and $[a,b]$ is any interval.

\paragraph{Proof.}We have for $\mu-$a.e. $t\in \Omega$, for $u_{1},u_{2} \in \left[ a,b \right]$, $u_{1}< u_{2}$ that
\begin{align}
p_{-}(t,u_{1}) \leq p_{+}(t,u_{1}) \leq p_{-}(t,u_{2}),
\end{align}
where 
\begin{align*}
p_{-}(t,u) = \lim \limits_{u \uparrow u_{o}} \frac{M(t,u)-M(t,u_{o})}{u-u_{o}}
\end{align*}
and 
\begin{align*}
p_{+}(t,u) = \lim \limits_{u \downarrow u_{o}} \frac{M(t,u)-M(t,u_{o})}{u-u_{o}},
\end{align*}
that is, for $\mu-$a.e. $t\in \Omega, p_{-}(t,u)$ is monotonic and hence it is continuous almost everywhere of $u$ on $\left[ a,b \right]$ (see theorem 4.19 [30]). For $\mu-$a.e. $t\in \Omega$, let $u_{1}$ be a continuity point of $p_{-}(t,u)$ and by taking the limit in $(2.1)$ as $u_{2} \rightarrow u_{1}$, we get by the Squeeze theorem for functions that for $\mu-$a.e. $t\in \Omega$, 
\begin{align*}
p_{-}(t,u_{1}) \leq p_{+}(t,u_{1}) \leq p_{-}(t,u_{1})
\end{align*}
which means that for $\mu-$a.e. $t\in \Omega,  p_{-}(t,u_{1}) = p_{+}(t,u_{1})$. Moreover, from $(2.1)$ we have by the fundamental theorem of calculus that for $\mu-$a.e. $t\in \Omega, \frac{\partial M(t,u)}{\partial u} = p(t,u) = p_{+}(t,u)$. Since for $\mu-$a.e. $t\in \Omega$, the function $M(t,u)$ is convex of $u$ on $\left[ a,b\right]$, then for $\mu-$a.e. $t\in \Omega, M(t,u)$ is absolutely continuous of $u$ on $\left[ a,b\right]$ (see lemma 1.3 [17]) and that for $\mu-$a.e. $t\in \Omega, M(t,u)$ is indefinite integral of its derivative $\frac{\partial M(t,u)}{\partial u}$ (see theorem 13.17 [32]), that is

\begin{align*}
M(t,u) = M(t,a) + \int_{a}^{u}\frac{\partial M(t,s)}{\partial s}ds = \int_{a}^{u}p(t,s)ds
\end{align*}
for $u \in [a,b].$
 																
\begin{flushright}
$\square$
\end{flushright}

\paragraph{Theorem 2.2.} The function $M:\Omega \times \mathbb{R} \to \mathbb{R} $ is Musielak $N$-function if and only if it can be written as follows: for $\mu-$a.e. $t\in \Omega$, 
\begin{align}
M(t,u) = \int_{0}^{\vert u \vert} p(t,s)ds,
\end{align}
where $p :\Omega \times [0,\infty ) \to \mathbb{R}$ is, for $\mu-$a.e. $t\in \Omega$, a non-decreasing, right-continuous function of $0$ on $\left[0,\infty\right)$ satisfies $p(t,0) = 0, p(t,u) > 0$ when $u>0$ and $\lim \limits_{u \to \infty} p(t,u) = \infty $ and for each $u \in \mathbb{R}$, $p(t,u)$ is a $\mu-$ measurable function of $t$ on $\Omega.$

\paragraph{Proof.} Given that $M(t,u)$ is a Musielak $N$-function, then for $\mu-$a.e. $t\in \Omega$, $M(t,u)$ is convex of $u$ on $\mathbb{R}$ and $M(t,0) = 0$. By theorem 2.1 that for $\mu-$a.e. $t\in \Omega,$

\begin{align*}
M(t,u) = \int_{0}^{\vert u \vert} p_{+}(t,s)ds \leq p_{+}(t,u)\int_{0}^{u}ds = up_{+}(t,u),
\end{align*}
i.e, 
\begin{align*} 
\frac{M(t,u)}{u} \leq p_{+}(t,u),
\end{align*}
where $p_{+}(t,u):\Omega \times \left[ 0,\infty \right)$ is, for $\mu-$a.e. $t\in \Omega$, a non-decreasing, right-continuous function of $0$ on $\left[ 0,\infty \right)$. From $(2)$ and $(4)$ of the definition 2.1, we have for $\mu-$a.e. $t\in \Omega$ that $p_{+}(t,u) > 0$ whenever $u>0$ and $\lim \limits_{u \to \infty} p_{+}(t,u)= \infty $ respectively. Moreover, for $\mu-$a.e. $t\in \Omega$ and for any $u > 0$ that
\begin{align*}
M(t,2u) &= \int_{0}^{2u} p_{+}(t,s)ds > \int_{u}^{2u} p_{+}(t,s)ds \\
		&> p_{+}(t,u) \int_{u}^{2u}ds = up_{+}(t,u),
\end{align*}
i.e, $$ p_{+}(t,u) \leq \frac{M(t,2u)}{u}.$$
\\
Again from $(3)$ of definition of Musielak $N$-function for $\mu-$a.e. $t\in \Omega$ that
\begin{align*}
p_{+}(t,0) &= \lim \limits_{u \to +0} p_{+}(t,u) \leq \lim \limits_{u \to +0} \frac{M(t,2u)}{u} \\
		   &= 2 \lim \limits_{u \to +0} \frac{M(t,2u)}{2u} = 0
\end{align*}
and for $\mu-$a.e. $t\in \Omega,$ 
\begin{align*}
 0 = \lim \limits_{u \to +0} \frac{M(t,u)}{u} \leq \lim \limits_{u \to +0} p_{+}(t,u)= p_{+}(t,0).  \end{align*}
Thus by the Squeeze theorem of functions that for $\mu-$a.e. $t\in \Omega, p_{+}(t,0) = 0$.
Since for $\mu-$a.e. $t\in \Omega$ that 

\begin{align*}
p_{+}(t,u) = \lim \limits_{u \downarrow u_{o}} \frac{M(t,u)-M(t,u_{o})}{u-u_{o}}
\end{align*}
and for each $u \in \mathbb{R}$ that $M(t,u)$ is a $\mu-$measurable function of $t$ on $\Omega$, then for each $u \in \mathbb{R}, p_{+}(t,u)$ is measurable function of $t$ on $\Omega$.

Now given that $M(t,u)$ satisfies $(2.3)$. So, for $\mu-$a.e. $t\in \Omega$, $M(t,u)$ is even and positive for any $u>0$. By theorem 2.1, for $\mu-$a.e. $t\in \Omega$, $M(t,u)$ is convex of $u$ on $\mathbb{R}$ and by $(2.3)$ for $\mu-$a.e $t\in \Omega$,
\begin{align*}
  0 < \frac{M(t,u)}{u} \leq p_{+}(t,u) < \frac{M(t,2u)}{u}
\end{align*}
for any $u>0$. Then
\begin{align*}
0 < \lim \limits_{u \to +0} \frac{M(t,u)}{u} \leq \lim \limits_{u \to 0} p_{+}(t,u) = 0.
\end{align*}
By the Squeeze theorem for functions that for $\mu-$a.e. $t\in \Omega$, 

\begin{align*}
 \lim \limits_{u \to 0} \frac{M(t,u)}{u} = 0
\end{align*} 
and
\begin{align*}
2 \lim \limits_{2u \to +\infty} \frac{M(t,2u)}{2u} >  \lim \limits_{u \to +\infty} \frac{M(t,u)}{u} = +\infty,
\end{align*} 
i.e, 
\begin{align*}
  \lim \limits_{u \to +\infty} \frac{M(t,u)}{u} = +\infty. 
\end{align*} 
Since for each $u \in \mathbb{R}$, $p_{+}(t,u)$ is measurable function of $t$ on $\Omega $ and for $\mu-$a.e. $t\in \Omega$, $M(t,u)$ satisfies $(2.3)$, then for each $u \in \mathbb{R}, M(t,u)$ is a $\mu-$measurable function of $t$ on $\Omega$. Therefore, $M:\Omega \times \mathbb{R} \to \mathbb{R} $ is Musielak $N$-function.

\begin{flushright}
$\square$
\end{flushright}

\paragraph{Theorem 2.3.} Any Musielak $N$-function $M:\Omega \times \mathbb{R} \to \mathbb{R}$ is continuous from the right of 0 on $\mathbb{R}$ for $\mu-$a.e. $t\in \Omega$.

\paragraph{Proof.} Given that $M:\Omega \times \mathbb{R} \to \mathbb{R}$ is Musielak $N$-function. By theorem 2.2, we have for $\mu-$a.e. $t\in \Omega$, 
\begin{eqnarray*}
0 \leq \lim \limits_{u \to +0} M(t,u) &=& \lim \limits_{u \to +0} \int_{0}^{\vert u \vert} p(t,s)ds = \int_{0}^{0} p(t,s)ds \\
                                &=& \lim \limits_{n \to \infty} \sum_{i = 1}^{n}p(t,u_{i})\Delta u 
                               = \lim \limits_{n \to \infty} \sum_{i = 1}^{n}p(t,u_{i})0 = 0
\end{eqnarray*}
where $u_{i} = 0 + i \Delta u$ and $\Delta u = \frac{0}{n}, n \in \mathbb{N}$. By the Squeeze theorem for functions that for $\mu-$a.e. $t\in \Omega, \lim \limits_{u \to +0} M(t,u) = 0.$

\begin{flushright}
$\square$
\end{flushright}

\section{Properties, calculus and basic approximation of Musielak $N$-functions} 
In this section, we are going to define the Musielak-Orlicz space generated by a Musielak $N$-function which is similar to the one that generated by a Musielak-Orlicz function [21] and investigate properties, calculus and basic approximation of Musielak $N$-functions and the Musielak-Orlicz spaces generated by them using the facts and results of the measure theory. 

\paragraph{Definition 3.1.} Let $\left( \Omega,\Sigma, \mu \right)$ be a measure space and $M$ be a Musielak $N$-function. The Musielak-Orlicz space $L_{M}\left( \Omega,\Sigma, \mu \right)$ generated by $M$ is defined by  
\begin{align*}
L_{M}\left( \Omega,\Sigma, \mu \right) = \left\lbrace f \in BS_{\Omega}: \int_{\Omega}M(t,\Vert f(t) \Vert_{BS})d\mu < \infty \right\rbrace
\end{align*}
where $BS_{\Omega}$ is the set of all $\mu-$measurable functions from $\Omega$ to $BS$ and $(BS,\Vert \cdot \Vert_{BS})$ is a Banach space, with Luxemburg norm 
\begin{align*}
\Vert f \Vert_{M} = \inf \left\lbrace  \lambda > 0:\int_{\Omega}M(t, \frac{\Vert f(t) \Vert_{BS}}{\lambda})d\mu \leq 1 \right\rbrace  
\end{align*}
 
\paragraph{Remark 3.1.} The Musielak-Orlicz space generated by a Musielak $N$-function $M$ on a measure space $(\Omega,\Sigma, \mu)$ is the Orlicz space generated by an $N$-function $\varphi$ on a measure space $(\Omega,\Sigma, \mu)$ whenever the Musielak $N$-function $M$ is the $N$-function $\varphi$. That is, $L_{M}(\Omega,\Sigma,\mu) = L_{\varphi}(\Omega,\Sigma,\mu)$ whenever for $\mu-$a.e. $t\in \Omega$ that $M(t,u) = \varphi(u)$ of $u$ on $\mathbb{R}$. 

\paragraph{Remark 3.2.}Every Musielak $N$-function $M$ is a Musielak-Orlicz function $MO$ with two additional conditions: for $\mu-$a.e. $t\in \Omega$,
\begin{eqnarray*}
\lim \limits_{u \to 0} \frac{M(t,u)}{u} = 0, \lim \limits_{u \to \infty} \frac{M(t,u)}{u} = \infty;
\end{eqnarray*}
so, the set of all Musielak $N$-functions 
\begin{eqnarray*}
\textbf{F}_{M} = \left\lbrace  M| M:\Omega \times \mathbb{R} \to \mathbb{R} \ \mathrm{is\ a\ Musielak} \ N-\mathrm{function} \right\rbrace 
\end{eqnarray*}
is contained in the set of all Musielak-Orlicz functions
\begin{eqnarray*}
\textbf{F}_{MO} = \left\lbrace  MO| MO:\Omega \times \left[ 0, \infty \right)  \to \left[ 0, \infty \right)  \ \mathrm{is\ a\ Musielak-Orlicz \ function} \right\rbrace. 
\end{eqnarray*}

\paragraph{Remark 3.3.}  The Musielak-Orlicz space generated by a Musielak $N$-function $M$ with the luxemburg norm $\left( L_{M}(\Omega,\Sigma,\mu),\Vert \cdot \Vert_{M}\right)$ is a Banach space, since the Musielak-Orlicz space generated by a Musielak-Orlicz function with the luxemburg norm is a Banach space (see theorem 7.7[21]) and by remark 3.2 the Musielak $N$-function is a Musielak-Orlicz function with two additional conditions.

\paragraph{Theorem 3.1.} Let $M_{1}$ and $M_{2}$ be two Musielak $N$-functions such that for $\mu-$a.e. $t\in \Omega, M_{2}(t,u) \leq r M_{1}(t,u)$ for some number $r>0$ and all $u \geq u_{0} > 0$. Then $L_{M_{1}}(\Omega,\Sigma,\mu) \subseteq L_{M_{2}}(\Omega,\Sigma,\mu).$ 

\paragraph{Proof.} Take $f \in L_{M_{1}}(\Omega,\Sigma,\mu),$ then for $\mu-$a.e. $t\in \Omega$,
\begin{align*}
\int_{\Omega}M_{1}(t,\Vert f(t) \Vert_{BS})d\mu < \infty.
\end{align*} 
From the assumption that for $\mu-$a.e. $t\in \Omega,$ there exist $r > 0$ and $u_{0} > 0$,  
\begin{align*}
M_{2}(t,u) \leq r M_{1}(t,u), \ u \geq u_{0},
\end{align*}
then for $\mu-$a.e. $t\in \Omega$, there exist $r > 0$ and $u_{0} > 0$ such that,
\begin{eqnarray*}
\int_{\Omega} M_{2}(t,\Vert f(t) \Vert_{BS}) d\mu \leq r \int_{\Omega} M_{1}(t,\Vert f(t) \Vert_{BS}) d\mu < \infty, \ f \neq 0
\end{eqnarray*}
i.e, $f \in L_{M_{2}}(\Omega,\Sigma,\mu).$
 
\begin{flushright}
$\square$
\end{flushright}

\paragraph{Corollary 3.1} Let $M_{1}$ and $M_{2}$ be two Musielak $N$-functions such that for $\mu-$a.e. $t\in \Omega,r_{1}M_{1}(t,u) \leq M_{2}(t,u) \leq r_{2} M_{1}(t,u)$ for some numbers $r_{1},r_{2}>0$ and all $u \geq u_{0} > 0$. Then $L_{M_{1}}(\Omega,\Sigma,\mu) = L_{M_{2}}(\Omega,\Sigma,\mu).$

\paragraph{Proof.} From the assumption there exist $r_{1},r_{2}>0$ and $u_{0} > 0$ such that for $\mu-$a.e. $t\in \Omega$, $M_{2}(t,u) \leq r_{2} M_{1}(t,u)$ and $M_{1}(t,u) \leq r_{3} M_{2}(t,u), r_{3} = 1/r_{1}$. From theorem 3.1 we have $L_{M_{1}}(\Omega,\Sigma,\mu) \subseteq L_{M_{2}}(\Omega,\Sigma,\mu)$ and $L_{M_{2}}(\Omega,\Sigma,\mu) \subseteq L_{M_{1}}(\Omega,\Sigma,\mu).$

\begin{flushright}
$\square$
\end{flushright}

\paragraph{Theorem 3.2.} If $M_{1}:\Omega \times \mathbb{R} \to \mathbb{R}$ is a Musielak $N$-function and $M_{2}:\Omega \times \mathbb{R} \to \mathbb{R}$ is a function such that for $\mu-$a.e. $t\in \Omega, M_{1}(t,u) = M_{2}(t,u)$ of $u$ on $\mathbb{R}$, then $M_{2}$ is a Musielak $N$-function

\paragraph{Proof.} It is given that $M_{1}$ is a Musielak $N$-function and for $\mu-$a.e. $t\in \Omega, M_{1}(t,u) = M_{2}(t,u)$ of $u$ on $\mathbb{R}$, then for $\mu-$a.e. $t\in \Omega, M_{2}(t,u) = M_{1}(t,u)$ is even convex function of $u$ on $\mathbb{R}$, for $\mu-$a.e. $t\in \Omega, M_{2}(t,u) = M_{1}(t,u) > 0$ for any $u > 0$ and for $\mu-$a.e. $t\in \Omega$,
\begin{eqnarray*}
\lim \limits_{u \to 0} \frac{M_{2}(t,u)}{u} = \lim \limits_{u \to 0} \frac{M_{1}(t,u)}{u} = 0
\end{eqnarray*}
and
\begin{eqnarray*}
\lim \limits_{u \to \infty} \frac{M_{2}(t,u)}{u} = \lim \limits_{u \to \infty} \frac{M_{1}(t,u)}{u} = \infty.
\end{eqnarray*}
Moreover, for each $u \in \mathbb{R}, M_{2}(t,u)$ is a $\mu-$measurable function of $t$ on $\Omega$. Then, $M_{2}$ is a Musielak $N$-function on $\Omega \times \mathbb{R}$ according to definition 2.1.

\begin{flushright} 
$\square$
\end{flushright}

\paragraph{Theorem 3.3.} If $\left( L_{M_{1}}(\Omega,\Sigma,\mu),\Vert \cdot \Vert_{M_{1}}\right)$ is a Musielak-Orlicz space generated by a Musielak $N$-function $M_{1}$ with the luxemburg norm and $M_{2}:\Omega \times \mathbb{R} \to \mathbb{R}$ is a function such that for $\mu-$a.e. $t\in \Omega, M_{1}(t,u) = M_{2}(t,u)$ of $u$ on $\mathbb{R}$, then for $\mu-$a.e. $t\in \Omega, \left( L_{M_{1}(t,u)}(\Omega,\Sigma,\mu),\Vert \cdot \Vert_{M_{1}(t,u)}\right) = \left( L_{M_{2}(t,u)}(\Omega,\Sigma,\mu),\Vert \cdot \Vert_{M_{2}(t,u)} \right)$ of $u$ on $\mathbb{R}$ and hence $\left( L_{M_{2}}(\Omega,\Sigma,\mu),\Vert \cdot \Vert_{M_{2}} \right)$ is a Musielak-Orlicz space generated by $M_{2}$ with the luxemburg norm.

\paragraph{Proof.} Since for $\mu-$a.e. $t\in \Omega, M_{1}(t,u) = M_{2}(t,u)$ of $u$ on $\mathbb{R}$, from theorem 3.2 we have that $M_{2}$ is a Musielak $N$-function and
\begin{eqnarray*}
L_{M_{2}(t,u)}\left( \Omega,\Sigma, \mu \right) &=& \left\lbrace f \in BS_{\Omega}: \int_{\Omega}M_{2}(t,\Vert f(t) \Vert_{BS})d\mu < \infty \right\rbrace \\
&=& \left\lbrace f \in BS_{\Omega}: \int_{\Omega}M_{1}(t,\Vert f(t) \Vert_{BS})d\mu < \infty \right\rbrace \\
&=& L_{M_{1}(t,u)}\left( \Omega,\Sigma, \mu \right) 
\end{eqnarray*} 
and for $f \in L_{M_{1}(t,u)}(\Omega,\Sigma,\mu) \left( = L_{M_{2}(t,u)}(\Omega,\Sigma,\mu) \right)$ that
\begin{eqnarray*}
\Vert f \Vert_{M_{1}(t,u)} &=& \inf \left\lbrace  \lambda > 0:\int_{\Omega}M_{1} \left( t, \frac{\Vert f(t) \Vert_{BS}}{\lambda}\right) d\mu \leq 1 \right\rbrace \\
 &=& \inf \left\lbrace  \lambda > 0:\int_{\Omega}M_{2} \left( t, \frac{\Vert f(t) \Vert_{BS}}{\lambda} \right) d\mu \leq 1 \right\rbrace \\
 &=& \Vert f \Vert_{M_{2}(t,u)};
\end{eqnarray*}
that is, for $\mu-$a.e. $t\in \Omega,\Vert f \Vert_{M_{2}(t,u)}$ can satisfies the norm's properties on $ L_{M_{2}(t,u)}(\Omega,\Sigma,\mu)$ of $u$ on $\mathbb{R}$. Then, for $\mu-$a.e. $t\in \Omega$, the equality hold in the assumption of $u$ on $\mathbb{R}$ and $ \left( L_{M_{2}}(\Omega,\Sigma,\mu),\Vert \cdot \Vert_{M_{2}}\right) $ is a Musielak-Orlicz space generated by $M_{2}$ with the luxemburg norm.

\begin{flushright}
$\square$
\end{flushright}  
 
\paragraph{Theorem 3.4.} If $\left\lbrace M_{n}: n \in \mathbb{N}\right\rbrace$ is a sequence of Musielak $N$-functions $M_{n}:\Omega \times \mathbb{R} \to \mathbb{R}$, then
\begin{align*}
\sup_{n \in \mathbb{N}} M_{n}, \inf_{n \in \mathbb{N}} M_{n}, \lim \limits_{n \to \infty} \sup  M_{n}, \lim \limits_{n \to \infty} \inf M_{n}
\end{align*}
are Musielak $N$-functions on $\Omega \times \mathbb{R}$ into $\mathbb{R}$.

\paragraph{Proof.} Since for all $n \in \mathbb{N}$ that $M_{n}$ is a Musielak $N$-function, we have for $\mu-$a.e. $t\in \Omega,  \sup_{n \in \mathbb{N}} M_{n}(t,u)$ and $\inf_{n \in \mathbb{N}} M_{n}(t,u)$ are even convex of $u$ on $\mathbb{R}$;
for $\mu-$a.e. $t\in \Omega, \sup_{n \in \mathbb{N}} M_{n}(t,u) > 0$ and $\inf_{n \in \mathbb{N}} M_{n}(t,u) > 0$ for any $u > 0$ and for $\mu-$a.e. $t\in \Omega$,
\begin{eqnarray*}
\lim \limits_{u \to 0} \frac{\sup_{n \in \mathbb{N}} M_{n}(t,u)}{u} = 0 &,& \lim \limits_{u \to + \infty} \frac{\sup_{n \in \mathbb{N}} M_{n}(t,u)}{u} = + \infty 
\end{eqnarray*}
and
\begin{eqnarray*}
\lim \limits_{u \to 0} \frac{\inf_{n \in \mathbb{N}} M_{n}(t,u)}{u} = 0  &,& \lim \limits_{u \to +\infty} \frac{\inf_{n \in \mathbb{N}} M_{n}(t,u)}{u} = +\infty. 
\end{eqnarray*} 
Moreover, that for each $u \in \mathbb{R}, \sup  M_{n}(t,u)$ and $\inf M_{n}(t,u)$ are $\mu-$measurable functions of $t$ on $\Omega$. So $\sup_{n \in \mathbb{N}} M_{n}$ and $\inf_{n \in \mathbb{N}} M_{n}$ are Musielak $N$-function on $\Omega \times \mathbb{R}$ . Also, for $\mu-$a.e. $t\in \Omega$, 
\begin{align*}
\lim \limits_{n \to \infty} \sup_{n \in \mathbb{N}}  M_{n}(t,u) = \inf_{n \in \mathbb{N}} \sup_{k \geq n}M_{k}(t,u), \\
\lim \limits_{n \to \infty} \inf_{n \in \mathbb{N}}  M_{n}(t,u) = \sup_{n \in \mathbb{N}} \inf_{k \geq n}M_{k}(t,u)
\end{align*}
of $u$ on $\mathbb{R}$, it follows that $\lim \limits_{n \to \infty} \sup  M_{n}$ and $\lim \limits_{n \to \infty} \inf M_{n}$ are Musielak $N$-functions on $\Omega \times \mathbb{R}$.

\begin{flushright}
$\square$
\end{flushright}
 
\paragraph{Theorem 3.5.} Let $\{M_{n}: n \in \mathbb{N}\}$ be a sequence of Musielak $N$-functions $M_{n}:\Omega \times \mathbb{R} \to \mathbb{R}$ such that for each $u \in \mathbb{R}, M_{n}(t,u) \leq M_{n+1}(t,u)$ on $\Omega$ for all $n \in \mathbb{N}$  . If $\left\lbrace \left( L_{M_{n}}(\Omega,\Sigma,\mu),\Vert \cdot \Vert_{M_{n}}\right): n \in \mathbb{N} \right\rbrace $ is a sequence of Musielak-Orlicz spaces generated by $\{ M_{n}: n \in \mathbb{N}\}$ with the luxemburg norm respectively, then $\sup_{n \in \mathbb{N}}( L_{M_{n}}(\Omega,\Sigma,\mu),\Vert \cdot \Vert_{M_{n}}) = ( L_{S}(\Omega,\Sigma,\mu),\Vert \cdot \Vert_{S})$, $\inf_{n \in \mathbb{N}}( L_{M_{n}}(\Omega,\Sigma,\mu),\Vert \cdot \Vert_{M_{n}}) = ( L_{I}(\Omega,\Sigma,\mu),\Vert \cdot \Vert_{I})$, $\lim \limits_{n \to \infty} \sup( L_{M_{n}}(\Omega,\Sigma,\mu),\Vert \cdot \Vert_{M_{n}}) = ( L_{LS}(\Omega,\Sigma,\mu),\Vert \cdot \Vert_{LS})$ and $\lim \limits_{n \to \infty} \inf(L_{M_{n}}(\Omega,\Sigma,\mu),\Vert \cdot \Vert_{M_{n}}) = ( L_{LI}(\Omega,\Sigma,\mu),\Vert \cdot \Vert_{LI})$ via $\{ M_{n}: n \in \mathbb{N}\}$, and hence $( L_{S}(\Omega,\Sigma,\mu),\Vert \cdot \Vert_{S})$, $( L_{I}(\Omega,\Sigma,\mu),\Vert \cdot \Vert_{I})$, $( L_{LS}(\Omega,\Sigma,\mu),\Vert \cdot \Vert_{LS})$ and $( L_{LI}(\Omega,\Sigma,\mu),\Vert \cdot \Vert_{LI})$ are Musielak-Orlicz spaces generated by 
\begin{align*}
S = \sup_{n \in \mathbb{N}} M_{n}, I = \inf_{n \in \mathbb{N}} M_{n}, LS = \lim \limits_{n \to \infty} \sup  M_{n} \ \mathrm{and} \ LI = \lim \limits_{n \to \infty} \inf M_{n}
\end{align*}
with the luxemburg norm respectively.

\paragraph{Proof.} We have from theorem 3.4 that $ S, I, LS $ and $LI$ are Musielak $N$-functions; so by the monotone convergence theorem that for $\mu-$a.e. $t\in \Omega$, 
\begin{eqnarray*}  
L_{M_{S}(t,u)}\left( \Omega,\Sigma, \mu \right) &=& \left\lbrace f \in BS_{\Omega}: \int_{\Omega}M_{S} \left( t,\Vert f(t) \Vert_{BS} \right) d\mu < \infty \right\rbrace \\
&=& \left\lbrace f \in BS_{\Omega}: \int_{\Omega} \sup_{n \in \mathbb{N}}M_{n} \left( t,\Vert f(t) \Vert_{BS} \right) d\mu < \infty \right\rbrace \\
&=& \left\lbrace f \in BS_{\Omega}: \sup_{n \in \mathbb{N}} \int_{\Omega}M_{n} \left( t,\Vert f(t) \Vert_{BS} \right) d\mu < \infty \right\rbrace \\
&=& \sup_{n \in \mathbb{N}} L_{M_{n}(t,u)}\left( \Omega,\Sigma, \mu \right), \\ 
L_{M_{I}(t,u)}\left( \Omega,\Sigma, \mu \right) &=& \left\lbrace f \in BS_{\Omega}: \int_{\Omega}M_{I} \left( t,\Vert f(t) \Vert_{BS} \right) d\mu < \infty \right\rbrace \\
&=& \left\lbrace f \in BS_{\Omega}: \int_{\Omega} \inf_{n \in \mathbb{N}}M_{n} \left( t,\Vert f(t) \Vert_{BS} \right) d\mu < \infty \right\rbrace \\
&=& \left\lbrace f \in BS_{\Omega}: \inf_{n \in \mathbb{N}} \int_{\Omega}M_{n} \left( t,\Vert f(t) \Vert_{BS} \right) d\mu < \infty \right\rbrace \\
&=& \inf_{n \in \mathbb{N}} L_{M_{n}(t,u)}\left( \Omega,\Sigma, \mu \right); 
\end{eqnarray*} 
and for $f,g \in L_{M_{S}(t,u)}(\Omega,\Sigma,\mu)(f,g \in L_{M_{I}(t,u)}(\Omega,\Sigma,\mu))$,
\begin{eqnarray*}
\Vert f \Vert_{M_{S}(t,u)} &=& \inf \left\lbrace  \lambda > 0:\int_{\Omega}M_{S} \left( t, \frac{\Vert f(t) \Vert_{BS}}{\lambda} \right) d\mu \leq 1 \right\rbrace \\
&=& \inf \left\lbrace  \lambda > 0:\int_{\Omega} \sup_{n \in \mathbb{N}}M_{n} \left( t, \frac{\Vert f(t) \Vert_{BS}}{\lambda} \right) d\mu \leq 1 \right\rbrace \\
&=& \inf \left\lbrace  \lambda > 0:\sup_{n \in \mathbb{N}} \int_{\Omega} M_{n} \left( t, \frac{\Vert f(t) \Vert_{BS}}{\lambda} \right) d\mu \leq 1 \right\rbrace \\
&=& \sup_{n \in \mathbb{N}} \Vert f \Vert_{M_{n}(t,u)}, \\
\Vert f \Vert_{M_{I}(t,u)} &=& \inf \left\lbrace  \lambda > 0:\int_{\Omega} M_{I} \left( t, \frac{\Vert f(t) \Vert_{BS}}{\lambda} \right) d\mu \leq 1 \right\rbrace \\
&=& \inf \left\lbrace  \lambda > 0:\int_{\Omega} \inf_{n \in \mathbb{N}} M_{n} \left( t, \frac{\Vert f(t) \Vert_{BS}}{\lambda} \right) d\mu \leq 1 \right\rbrace \\
&=& \inf \left\lbrace  \lambda > 0:\inf_{n \in \mathbb{N}} \int_{\Omega} M_{n} \left( t, \frac{\Vert f(t) \Vert_{BS}}{\lambda} \right) d\mu \leq 1 \right\rbrace \\
&=& \inf_{n \in \mathbb{N}} \Vert f \Vert_{M_{n}(t,u)};
\end{eqnarray*}
so for $\mu-$a.e. $t\in \Omega$, that $\Vert f \Vert_{M_{S}(t,u)} = \sup_{n \in \mathbb{N}} \Vert f \Vert_{M_{n}(t,u)} \geq 0$ and $\Vert f \Vert_{M_{I}(t,u)} = \inf_{n \in \mathbb{N}} \Vert f \Vert_{M_{n}(t,u)} \geq 0$ of $u$ on $\mathbb{R}$. For $\mu-$a.e. $t\in \Omega$, $\Vert f \Vert_{M_{S}(t,u)} = 0 \Rightarrow \sup_{n \in \mathbb{N}} \Vert f \Vert_{M_{n}(t,u)} = 0$ if and only if $f(t) = 0$ and $\Vert f \Vert_{M_{I}(t,u)} = 0 \Rightarrow \inf_{n \in \mathbb{N}} \Vert f \Vert_{M_{n}(t,u)} = 0$ if and only if $f(t) = 0$ for an arbitrary $\lambda$. For $\mu-$a.e. $t\in \Omega$, for any scalar $\alpha, \Vert \alpha f \Vert_{M_{S}(t,u)} = \sup_{n \in \mathbb{N}} \Vert \alpha f \Vert_{M_{n}(t,u)} = |\alpha| \sup_{n \in \mathbb{N}} \Vert f \Vert_{M_{n}(t,u)} = |\alpha| \Vert f \Vert_{M_{S}(t,u)}$ and $\Vert \alpha f \Vert_{M_{I}(t,u)} = \inf_{n \in \mathbb{N}}\Vert \alpha f\Vert_{M_{n}(t,u)} = |\alpha| \inf_{n \in \mathbb{N}} \Vert f \Vert_{M_{n}(t,u)} = |\alpha| \Vert f \Vert_{M_{I}(t,u)}$ of $u$ on $\mathbb{R}$. And since for all $n \in \mathbb{N}$, for $\mu-$a.e. $t\in \Omega,M_{n}(t,u)$ are convex functions of $u$ on $\mathbb{R}$, then for $\mu-$a.e. $t\in \Omega$ that
\begin{eqnarray*}
\int_{\Omega}M_{S}\left( t, \frac{\Vert f(t)+g(t)\Vert_{BS}}{\Vert f(t)\Vert_{M_{S}} + \Vert g(t)\Vert_{M_{S}}} \right) d\mu &=& \int_{\Omega} \sup_{n \in \mathbb{N}} M_{n} \left( t, \frac{\Vert f(t)+g(t)\Vert_{BS}}{ \Vert f(t)\Vert_{M_{S}} + \Vert g(t)\Vert_{M_{S}} } \right) d\mu \\
&=& \sup_{n \in \mathbb{N}} \int_{\Omega} M_{n} \left( t, \frac{\Vert f(t)+g(t)\Vert_{BS}}{ \Vert f(t)\Vert_{M_{S}} + \Vert g(t)\Vert_{M_{S}} } \right) d\mu \\
& \leq & \frac{\Vert f(t) \Vert_{M_{S}}}{\Vert f(t)\Vert_{M_{S}} + \Vert g(t)\Vert_{M_{S}}}\sup_{n \in \mathbb{N}} \int_{\Omega} M_{n} \left( t, \frac{\Vert f(t)\Vert_{BS}}{\Vert f(t)\Vert_{M_{S}} + \Vert g(t)\Vert_{M_{S}}} \right) d\mu \\
&+& \frac{\Vert g(t) \Vert_{M_{S}}}{\Vert f(t)\Vert_{M_{S}} + \Vert g(t)\Vert_{M_{S}}} \sup_{n \in \mathbb{N}} \int_{\Omega}  M_{n} \left( t, \frac{\Vert g(t)\Vert_{BS}}{\Vert f(t)\Vert_{M_{S}} + \Vert g(t)\Vert_{M_{S}}} \right) d\mu \\
& \leq & \frac{\Vert f(t) \Vert_{M_{S}}}{\Vert f(t)\Vert_{M_{S}} + \Vert g(t)\Vert_{M_{S}}}\sup_{n \in \mathbb{N}} \int_{\Omega} M_{n} \left( t, \frac{\Vert f(t)\Vert_{BS}}{\Vert f(t)\Vert_{M_{S}}} \right) d\mu \\
&+& \frac{\Vert g(t) \Vert_{M_{S}}}{\Vert f(t)\Vert_{M_{S}} + \Vert g(t)\Vert_{M_{S}}} \sup_{n \in \mathbb{N}} \int_{\Omega}  M_{n} \left( t, \frac{\Vert g(t)\Vert_{BS}}{\Vert g(t)\Vert_{M_{S}}} \right) d\mu \\
& \leq & 1,
\end{eqnarray*}
because $\int_{\Omega}  M_{n} \left( t, \frac{\Vert f(t)\Vert_{BS}}{\Vert f(t)\Vert_{M_{S}}} \right) d\mu \leq 1$ and $\int_{\Omega}  M_{n} \left( t, \frac{\Vert g(t)\Vert_{BS}}{\Vert g(t)\Vert_{M_{S}}} \right) d\mu \leq 1$, so for $\mu-$a.e. $t\in \Omega, \Vert f(t)+g(t)\Vert_{M_{S}(t,u)} \leq \Vert f(t)\Vert_{M_{S}(t,u)} + \Vert g(t)\Vert_{M_{S}(t,u)}$ of $u$ on $\mathbb{R}$;
and 
\begin{eqnarray*}
\int_{\Omega}M_{I}\left( t, \frac{\Vert f(t)+g(t)\Vert_{BS}}{\Vert f(t)\Vert_{M_{I}} + \Vert g(t)\Vert_{M_{I}}} \right) d\mu &=& \int_{\Omega} \inf_{n \in \mathbb{N}} M_{n} \left( t, \frac{\Vert f(t)+g(t)\Vert_{BS}}{ \Vert f(t)\Vert_{M_{I}} + \Vert g(t)\Vert_{M_{I}} } \right) d\mu \\
&=& \inf_{n \in \mathbb{N}} \int_{\Omega} M_{n} \left( t, \frac{\Vert f(t)+g(t)\Vert_{BS}}{ \Vert f(t)\Vert_{M_{I}} + \Vert g(t)\Vert_{M_{I}} } \right) d\mu \\
& \leq & \frac{\Vert f(t) \Vert_{M_{I}}}{\Vert f(t)\Vert_{M_{I}} + \Vert g(t)\Vert_{M_{I}}}\inf_{n \in \mathbb{N}} \int_{\Omega} M_{n} \left( t, \frac{\Vert f(t)\Vert_{BS}}{\Vert f(t)\Vert_{M_{I}} + \Vert g(t)\Vert_{M_{I}}} \right) d\mu \\
&+& \frac{\Vert g(t) \Vert_{M_{I}}}{\Vert f(t)\Vert_{M_{I}} + \Vert g(t)\Vert_{M_{I}}} \inf_{n \in \mathbb{N}} \int_{\Omega}  M_{n} \left( t, \frac{\Vert g(t)\Vert_{BS}}{\Vert f(t)\Vert_{M_{I}} + \Vert g(t)\Vert_{M_{I}}} \right) d\mu \\
& \leq & \frac{\Vert f(t) \Vert_{M_{I}}}{\Vert f(t)\Vert_{M_{I}} + \Vert g(t)\Vert_{M_{I}}}\inf_{n \in \mathbb{N}} \int_{\Omega} M_{n} \left( t, \frac{\Vert f(t)\Vert_{BS}}{\Vert f(t)\Vert_{M_{I}}} \right) d\mu \\
&+& \frac{\Vert g(t) \Vert_{M_{I}}}{\Vert f(t)\Vert_{M_{I}} + \Vert g(t)\Vert_{M_{I}}} \inf_{n \in \mathbb{N}} \int_{\Omega}  M_{n} \left( t, \frac{\Vert g(t)\Vert_{BS}}{\Vert g(t)\Vert_{M_{I}}} \right) d\mu \\
& \leq & 1,
\end{eqnarray*} 
because $\int_{\Omega}  M_{n} \left( t, \frac{\Vert f(t)\Vert_{BS}}{\Vert f(t)\Vert_{M_{I}}} \right) d\mu \leq 1$ and $\int_{\Omega}  M_{n} \left( t, \frac{\Vert g(t)\Vert_{BS}}{\Vert g(t)\Vert_{M_{I}}} \right) d\mu \leq 1$, so for $\mu-$a.e. $t\in \Omega, \Vert f(t)+g(t)\Vert_{M_{I}(t,u)} \leq \Vert f(t)\Vert_{M_{I}(t,u)} + \Vert g(t)\Vert_{M_{I}(t,u)}$ of $u$ on $\mathbb{R}$; that is, for $\mu-$a.e. $t\in \Omega, \Vert f(t)\Vert_{M_{S}(t,u)}$ and $\Vert f(t)\Vert_{M_{I}(t,u)}$ are norms on $L_{M_{S}(t,u)}(\Omega,\Sigma,\mu)$ and $L_{M_{I}(t,u)}(\Omega,\Sigma,\mu)$ of $u$ on $\mathbb{R}$ respectively. So, for $\mu-$a.e. $t\in \Omega$, the equality hold in the assumption of $u$ on $\mathbb{R}$ and $( L_{M_{S}}(\Omega,\Sigma,\mu),\Vert \cdot \Vert_{M_{S}})$ and $( L_{M_{I}}(\Omega,\Sigma,\mu),\Vert \cdot \Vert_{M_{I}}) $ are Musielak-Orlicz spaces generated by $M_{S}$ and $M_{I}$ with the luxemburg norm respectively. 

Now, again by the monotone convergence theorem that for $\mu-$a.e. $t\in \Omega$, 
\begin{eqnarray*}
L_{M_{LS}(t,u)}\left( \Omega,\Sigma, \mu \right) &=& \left\lbrace f \in BS_{\Omega}: \int_{\Omega}M_{LS}(t,\Vert f(t) \Vert_{BS})d\mu < \infty \right\rbrace \\
&=& \left\lbrace f \in BS_{\Omega}: \int_{\Omega} \lim \limits_{n \to \infty} \sup_{n \in \mathbb{N}}M_{n}(t,\Vert f(t) \Vert_{BS})d\mu < \infty \right\rbrace \\
&=& \left\lbrace f \in BS_{\Omega}: \lim \limits_{n \to \infty} \sup_{n \in \mathbb{N}} \int_{\Omega} M_{n}(t,\Vert f(t) \Vert_{BS})d\mu < \infty \right\rbrace \\
&=& \lim \limits_{n \to \infty} \sup_{n \in \mathbb{N}} L_{M_{k}(t,u)}\left( \Omega,\Sigma, \mu \right), \\
L_{M_{LI}(t,u)}\left( \Omega,\Sigma, \mu \right) &=& \left\lbrace f \in BS_{\Omega}: \int_{\Omega}M_{LI}(t,\Vert f(t) \Vert_{BS})d\mu < \infty \right\rbrace \\
&=& \left\lbrace f \in BS_{\Omega}: \int_{\Omega} \lim \limits_{n \to \infty} \inf_{n \in \mathbb{N}}M_{n}(t,\Vert f(t) \Vert_{BS})d\mu < \infty \right\rbrace \\
&=& \left\lbrace f \in BS_{\Omega}: \lim \limits_{n \to \infty} \inf_{n \in \mathbb{N}} \int_{\Omega} M_{n}(t,\Vert f(t) \Vert_{BS})d\mu < \infty \right\rbrace \\
&=& \lim \limits_{n \to \infty} \inf_{n \in \mathbb{N}}L_{M_{k}(t,u)}\left( \Omega,\Sigma, \mu \right); 
\end{eqnarray*} 
and for $f \in L_{M_{LS}(t,u)}(\Omega,\Sigma,\mu)$ $ \left( f \in L_{M_{LI}(t,u)}(\Omega,\Sigma,\mu)\right) $
\begin{eqnarray*}
\Vert f \Vert_{M_{LS}(t,u)} &=& \inf \left\lbrace  \lambda > 0:\int_{\Omega}M_{LS} \left( t, \frac{\Vert f(t) \Vert_{BS}}{\lambda}\right) d\mu \leq 1 \right\rbrace \\
&=& \inf \left\lbrace  \lambda > 0:\int_{\Omega} \lim \limits_{n \to \infty} \sup_{n \in \mathbb{N}} M_{n} \left( t, \frac{\Vert f(t) \Vert_{BS}}{\lambda} \right) d\mu \leq 1 \right\rbrace \\
&=& \inf \left\lbrace  \lambda > 0:\lim \limits_{n \to \infty} \sup_{n \in \mathbb{N}} \int_{\Omega} M_{n} \left( t, \frac{\Vert f(t) \Vert_{BS}}{\lambda} \right) d\mu \leq 1 \right\rbrace \\
&=& \lim \limits_{n \to \infty} \sup_{n \in \mathbb{N}} \Vert f \Vert_{M_{n}(t,u)}, \\
\Vert f \Vert_{M_{LI}(t,u)} &=& \inf \left\lbrace  \lambda > 0:\int_{\Omega}M_{LI} \left( t, \frac{\Vert f(t) \Vert_{BS}}{\lambda} \right) d\mu \leq 1 \right\rbrace \\
&=& \inf \left\lbrace  \lambda > 0:\int_{\Omega} \lim \limits_{n \to \infty} \inf_{n \in \mathbb{N}}M_{n} \left( t, \frac{\Vert f(t) \Vert_{BS}}{\lambda}\right) d\mu \leq 1 \right\rbrace \\
&=& \inf \left\lbrace  \lambda > 0:\lim \limits_{n \to \infty} \inf_{n \in \mathbb{N}} \int_{\Omega} M_{n} \left( t, \frac{\Vert f(t) \Vert_{BS}}{\lambda}\right) d\mu \leq 1 \right\rbrace \\
&=& \lim \limits_{n \to \infty} \inf_{n \in \mathbb{N}} \Vert f \Vert_{M_{n}(t,u)};
\end{eqnarray*}
that is, for $\mu-$a.e. $t\in \Omega,$ $\Vert f \Vert_{M_{LS}(t,u)}$ and $\Vert f \Vert_{M_{LI}(t,u)}$ can satisfy the norm's properties on $\left( L_{M_{LS}(t,u)}(\Omega,\Sigma,\mu) \right)$ and $\left( L_{M_{LI}(t,u)}(\Omega,\Sigma,\mu) \right)$ of $u$ on $\mathbb{R}$ respectively. Therefore, for $\mu-$a.e. $t\in \Omega$, we get the equality in the assumption of $u$ on $\mathbb{R}$ and $\left( L_{M_{LS}}(\Omega,\Sigma,\mu),\Vert \cdot \Vert_{M_{LS}}\right)$ and $\left( L_{M_{LI}}(\Omega,\Sigma,\mu),\Vert \cdot \Vert_{M_{LI}}\right)$ are Musielak-Orlicz spaces generated by $M_{LS}$ and $M_{LI}$ with the luxemburg norm respectively.

\begin{flushright}
$\square$
\end{flushright}

\paragraph{Theorem 3.6.} If $\left\lbrace M_{n}: n \in \mathbb{N} \right\rbrace$ is a sequence of Musielak $N$-functions that satisfy the $\bigtriangleup_{2}$-condition $M_{n}:\Omega \times \mathbb{R} \to \mathbb{R}$, and $M_{n} \to M, M:\Omega \times \mathbb{R} \to \mathbb{R}$ pointwisely as $n \to \infty$, then $M$ is a Musielak $N$-function and satisfy the $\bigtriangleup_{2}$-condition.

\paragraph{Proof.} Since the convergence of $M_{n}$ to $M$ is pointwisely, then for $\mu-$a.e. $t\in \Omega$ that
\begin{eqnarray*}
M(t,u) = \lim \limits_{n \to \infty} \inf_{n \in \mathbb{N}}M_{n}(t,u) = \lim \limits_{n \to \infty} \sup_{n \in \mathbb{N}}M_{n}(t,u)
\end{eqnarray*}
of $u$ on $\mathbb{R}$, and the $\bigtriangleup_{2}$-condition is clear to satisfy.
 
\begin{flushright}
$\square$
\end{flushright}

\paragraph{Theorem 3.7.} Let $\{M_{n}: n \in \mathbb{N}\}$ be a sequence of Musielak $N$-functions $M_{n}:\Omega \times \mathbb{R} \to \mathbb{R}$ such that for each $u \in \mathbb{R}, M_{n}(t,u) \to M(t,u)$ on $\Omega$ as $n \to \infty$ and $\vert M_{n}(t,u) \vert \leq G(t,u)$, where for each $u \in \mathbb{R}, G$ is absolutely integrable on $\Omega$. If $\left\lbrace \left( L_{M_{n}}(\Omega,\Sigma,\mu),\Vert \cdot \Vert_{M_{n}} \right): n \in \mathbb{N} \right\rbrace $ is a sequence of Musielak-Orlicz spaces generated by $ \left\lbrace M_{n}: n \in \mathbb{N}\right\rbrace $ with the luxemburg norm, then $ \left( L_{M_{n}}(\Omega,\Sigma,\mu),\Vert \cdot \Vert_{M_{n}} \right) \to \left( L_{M}(\Omega,\Sigma,\mu),\Vert \cdot \Vert_{M} \right)$ via $\left\lbrace M_{n}: n \in \mathbb{N}\right\rbrace $ as $n \to \infty$ and hence $\left( L_{M}(\Omega,\Sigma,\mu),\Vert \cdot \Vert_{M} \right)$ is a Musielak-Orlicz space generated by $M$ with the luxemburg norm.

\paragraph{Proof.} From the assumptions and from theorem 3.6 that $M$ is a Musielak $N$-function and so by the Lebesgue's dominated convergence theorem that for $\mu-$a.e. $t\in \Omega$,
\begin{eqnarray*}
L_{M(t,u)}\left( \Omega,\Sigma, \mu \right) &=& \left\lbrace f \in BS_{\Omega}: \int_{\Omega}M(t,\Vert f(t) \Vert_{BS})d\mu < \infty \right\rbrace \\
&=& \left\lbrace f \in BS_{\Omega}: \int_{\Omega} \lim \limits_{n \to \infty} M_{n}(t,\Vert f(t) \Vert_{BS})d\mu < \infty \right\rbrace \\
&=& \left\lbrace f \in BS_{\Omega}: \lim \limits_{n \to \infty} \int_{\Omega} M_{n}(t,\Vert f(t) \Vert_{BS})d\mu < \infty \right\rbrace \\
&=& \lim \limits_{n \to \infty} L_{M_{n}(t,u)}\left( \Omega,\Sigma, \mu \right);
\end{eqnarray*}
and for $f,g \in L_{M(t,u)}(\Omega,\Sigma,\mu)$,
\begin{eqnarray*}
\Vert f \Vert_{M(t,u)} &=& \inf \left\lbrace  \lambda > 0:\int_{\Omega}M \left( t, \frac{\Vert f(t) \Vert_{BS}}{\lambda} \right) d\mu \leq 1 \right\rbrace \\
&=& \inf \left\lbrace  \lambda > 0:\int_{\Omega} \lim \limits_{n \to \infty} M_{n} \left( t, \frac{\Vert f(t) \Vert_{BS}}{\lambda} \right) d\mu \leq 1 \right\rbrace \\
&=& \inf \left\lbrace  \lambda > 0:\lim \limits_{n \to \infty} \int_{\Omega} M_{n} \left( t, \frac{\Vert f(t) \Vert_{BS}}{\lambda} \right) d\mu \leq 1 \right\rbrace \\
&=& \lim \limits_{n \to \infty} \Vert f \Vert_{M_{n}(t,u)} ;
\end{eqnarray*}
that is, for $\mu-$a.e. $t\in \Omega$, that $\Vert f \Vert_{M(t,u)} = \lim \limits_{n \to \infty} \Vert f \Vert_{M_{n}(t,u)} \geq 0$ of $u$ on $\mathbb{R}$. For $\mu-$a.e. $t\in \Omega$, $\Vert f \Vert_{M(t,u)} = 0 \Rightarrow\lim \limits_{n \to \infty} \Vert f \Vert_{M_{n}(t,u)} = 0$ if and only if $f(t) = 0$ for an arbitrary $\lambda$. For $\mu-$a.e. $t\in \Omega$, for any scalar $\alpha, \Vert \alpha f \Vert_{M} = \lim \limits_{n \to \infty} \Vert \alpha f \Vert_{M_{n}} = \vert \alpha \vert \lim \limits_{n \to \infty} \Vert  f \Vert_{M_{n}} = \vert \alpha \vert \Vert  f \Vert_{M}$ of $u$ on $\mathbb{R}$. And since for all $n \in \mathbb{N}$, for $\mu-$a.e. $t\in \Omega, M_{n}(t,u)$ are convex of $u$ on $\mathbb{R}$, then for $f,g \in L_{M}\left( \Omega,\Sigma, \mu \right)$ we have that 
\begin{eqnarray*}
\int_{\Omega}M \left( t, \frac{\Vert f(t)+g(t)\Vert_{BS}}{\Vert f(t)\Vert_{M} + \Vert g(t)\Vert_{M}} \right) d\mu &=& \int_{\Omega} \lim \limits_{n \to \infty} M_{n} \left( t, \frac{\Vert f(t)+g(t)\Vert_{BS}}{ \Vert f(t)\Vert_{M} + \Vert g(t)\Vert_{M}} \right) d\mu \\
&=& \lim \limits_{n \to \infty} \int_{\Omega} M_{n} \left( t, \frac{\Vert f(t)+g(t)\Vert_{BS}}{ \Vert f(t)\Vert_{M} + \Vert g(t)\Vert_{M}} \right) d\mu \\
& \leq & \frac{\Vert f(t) \Vert_{M}}{\Vert f(t)\Vert_{M} + \Vert g(t)\Vert_{M}} \lim \limits_{n \to \infty} \int_{\Omega} M_{n} \left( t, \frac{\Vert f(t)\Vert_{BS}}{\Vert f(t)\Vert_{M} + \Vert g(t)\Vert_{M}} \right) d\mu \\
&+& \frac{\Vert g(t) \Vert_{M}}{\Vert f(t)\Vert_{M} + \Vert g(t)\Vert_{M}} \lim \limits_{n \to \infty} \int_{\Omega} M_{n} \left( t, \frac{\Vert g(t)\Vert_{BS}}{\Vert f(t)\Vert_{M} + \Vert g(t)\Vert_{M}} \right) d\mu \\
& \leq & \frac{\Vert f(t) \Vert_{M}}{\Vert f(t)\Vert_{M} + \Vert g(t)\Vert_{M}} \lim \limits_{n \to \infty} \int_{\Omega} M_{n} \left( t, \frac{\Vert f(t)\Vert_{BS}}{\Vert f(t)\Vert_{M}} \right) d\mu \\
&+& \frac{\Vert g(t) \Vert_{M}}{\Vert f(t)\Vert_{M} + \Vert g(t)\Vert_{M}} \lim \limits_{n \to \infty} \int_{\Omega} M_{n} \left( t, \frac{\Vert g(t)\Vert_{BS}}{\Vert g(t)\Vert_{M}} \right) d\mu \\
& \leq & 1,
\end{eqnarray*} 
because $\int_{\Omega} M_{n} \left( t, \frac{\Vert f(t)\Vert_{BS}}{\Vert f(t)\Vert_{M}} \right) \leq 1$  and $\int_{\Omega} M_{n} \left( t, \frac{\Vert g(t)\Vert_{BS}}{ \Vert g(t)\Vert_{M}} \right) \leq 1$, so for $\mu-$a.e. $t\in \Omega, \Vert f(t) + g(t) \Vert_{M(t,u)} \leq \Vert f(t) \Vert_{M(t,u)} + \Vert g(t) \Vert_{M(t,u)}$ of $u$ on $\mathbb{R}$; that is, for $\mu-$a.e. $t\in \Omega, \Vert f(t) \Vert_{M(t,u)}$ is norm on $L_{M(t,u)}\left( \Omega,\Sigma, \mu \right)$ of $u$ on $\mathbb{R}$. So, for $\mu-$a.e. $t\in \Omega$ the convergence holds in the assumption of $u$ on $\mathbb{R}$ and $\left( L_{M}(\Omega,\Sigma,\mu),\Vert \cdot \Vert_{M} \right)$ is a Musielak-Orlicz spaces generated by $M$ with the luxemburg norm.

\begin{flushright}
$\square$
\end{flushright}

\paragraph{Corollary 3.2.} Under theorem 3.7's assumptions with $\vert M_{n}(t,u) \vert \leq M(t,u)$ and $M_{n}$ satisfies the $\bigtriangleup_{2}$-condition for all $n \in \mathbb{N}$, if $f \in \left( L_{M}(\Omega,\Sigma,\mu),\Vert \cdot \Vert_{M} \right)$ then there exists a sequence of functions $\{f_{n}:n \in \mathbb{N}\}$  such that $f_{n} \in \left( L_{M_{n}}(\Omega,\Sigma,\mu),\Vert \cdot \Vert_{M_{n}} \right)$ for all $n \in \mathbb{N}$ and $f_{n} \underset{n \to \infty} \longrightarrow f$ under the luxemburg norm $\Vert \cdot \Vert_{M}$.

\paragraph{Proof.} From above assumptions and according to theorem 3.1 we have $L_{M}(\Omega,\Sigma,\mu) \subseteq L_{M_{n}}(\Omega,\Sigma,\mu)$ for all $n \in \mathbb{N}$; that is $ L_{M}(\Omega,\Sigma,\mu) = \bigcap_{n=1}^{\infty} L_{M_{n}}(\Omega,\Sigma,\mu)$. So, if $f \in L_{M}(\Omega,\Sigma,\mu) $, then $f \in L_{M_{n}}(\Omega,\Sigma,\mu)$ for all $n \in \mathbb{N}$. Fix $n_{0} \in \mathbb{N}$, since $\left( L_{M_{n_{0}}}(\Omega,\Sigma,\mu),\Vert \cdot \Vert_{M_{n_{0}}} \right)$ is separable, because $M_{n_{0}}$ satisfies the $\bigtriangleup_{2}$-condition[18], there exists $f_{n_{0}} \in \left( L_{M_{n_{0}}}(\Omega,\Sigma,\mu),\Vert \cdot \Vert_{M_{n_{0}}} \right)$ such that $ \Vert f_{n_{0}} - f \Vert_{M_{n_{0}}} < \frac{1}{n_{0}} $. Then,
\begin{eqnarray*}
0 \leq \Vert f_{n_{0}} - f \Vert_{M} = \lim \limits_{n \to \infty} \Vert f_{n_{0}} - f \Vert_{M_{n}} < \frac{1}{n_{0}} \lim \limits_{n|n_{0} \to \infty} \Vert f_{n_{0}} - f \Vert_{M_{n}}.
\end{eqnarray*} 
Letting $n_{0} \to \infty$, we get $\lim \limits_{n_{0} \to \infty} \Vert f_{n_{0}} - f \Vert_{M} = 0 $ by the Squeeze theorem of functions.

\begin{flushright}
$\square$
\end{flushright}

\paragraph{Theorem 3.8.} If $M_{i}:\Omega \times \mathbb{R} \to \mathbb{R}, i=1,2$ are Musielak $N$-functions, then for r $\in \mathbb{R_{+}}$ that $rM_{1}$ and $M_{1}+M_{2}$ are Musielak $N$-functions.

\paragraph{Proof.} Since $M_{i}, i=1,2$ are Musielak $N$-functions, then for $r \in \mathbb{R_{+}}$, for $\mu-$a.e. $t\in \Omega, rM_{1}(t,u)$ and $(M_{1}+M_{2})(t,u)$ are even convex of $u$ on $\mathbb{R}$; for $\mu-$a.e. $t\in \Omega, rM_{1}(t,u) > 0$ and $(M_{1}+M_{2})(t,u) > 0$  for any $u > 0$ and for $\mu-$a.e. $t\in \Omega$,

\begin{eqnarray*}
\lim \limits_{u \to 0} \frac{rM_{1}(t,u)}{u} = 0 ,& \lim \limits_{u \to +\infty } \frac{rM_{1}(t,u)}{u} = +\infty  \\ 
\lim \limits_{u \to 0} \frac{(M_{1}+M_{2})(t,u)}{u} = 0 ,& \lim \limits_{u \to +\infty} \frac{(M_{1}+M_{2})(t,u)}{u} = +\infty 
\end{eqnarray*}
Moreover, that for each $u \in \mathbb{R}, rM_{1}(t,u)$ and $(M_{1}+M_{2})(t,u)$ are $\mu-$measurable functions of $t$ on $\Omega$. So, $\Omega,rM_{1}$ and $M_{1}+M_{2}$ are Musielak $N$-functions on $\Omega \times \mathbb{R}$. 

\begin{flushright}
$\square$
\end{flushright}

\paragraph{Theorem 3.9.}If $\left( L_{M_{i}}(\Omega,\Sigma,\mu),\Vert \cdot \Vert_{M_{i}} \right): i = 1,2$ are Musielak-Orlicz spaces generated by Musielak $N$-functions $M_{i}: i = 1,2$  with the luxemburg norm respectively, then $ \left( L_{rM_{1}}(\Omega,\Sigma,\mu),\Vert \cdot \Vert_{rM_{1}} \right), r \geq 1$ and $\left( L_{M_{1}+M_{2}}(\Omega,\Sigma,\mu),\Vert \cdot \Vert_{M_{1}+M_{2}} \right)$ are Musielak-Orlicz spaces generated by Musielak $N$-functions $rM_{1}$ and $M_{1}+M_{2}$ with the luxemburg norm respectively.

\paragraph{Proof.} We have from theorem 3.8 for $r \geq 1$ that $rM_{1}$ and $M_{1}+M_{2}$ are Musielak $N$-functions, then for $\mu-$a.e. $t\in \Omega$, 
 
\begin{eqnarray*}
L_{rM_{1}(t,u)}(\Omega,\Sigma,\mu) &=& \left\lbrace f \in BS_{\Omega}: \int_{\Omega}rM_{1} \left( t,\Vert f(t) \Vert_{BS} \right) d\mu < \infty \right\rbrace , \\
&=& L_{M_{1}(t,u)}(\Omega,\Sigma,\mu), \\
L_{(M_{1}+M_{2})(t,u)}(\Omega,\Sigma,\mu) &=& \left\lbrace f \in BS_{\Omega}: \int_{\Omega}(M_{1}+M_{2}) \left( t,\Vert f(t) \Vert_{BS} \right) d\mu < \infty \right\rbrace;
\end{eqnarray*}
and for $f,g \in L_{rM_{1}(t,u)}(\Omega,\Sigma,\mu) \left( f,g \in L_{(M_{1}+M_{2})(t,u)}(\Omega,\Sigma,\mu) \right)$,

\begin{eqnarray*}
\Vert f \Vert_{rM_{1}(t,u)} &=& \inf \left\lbrace  \lambda > 0:\int_{\Omega}rM_{1} \left(t, \frac{\Vert f(t) \Vert_{BS}}{\lambda}\right)d\mu \leq 1 \right\rbrace \\
&=& \Vert f \Vert_{M_{1}(t,u)},\\
\Vert f \Vert_{(M_{1}+M_{2})(t,u)} &=& \inf \left\lbrace  \lambda > 0:\int_{\Omega}(M_{1}+M_{2})\left( t, \frac{\Vert f(t) \Vert_{BS}}{\lambda} \right) d\mu \leq 1 \right\rbrace.
\end{eqnarray*}
It is clear that $\Vert f \Vert_{rM_{1}(t,u)}$ is a norm on $L_{rM_{1}(t,u)}(\Omega,\Sigma,\mu)$. Now, for $\mu-$a.e. $t\in \Omega, \Vert f \Vert_{(M_{1}+M_{2})(t,u)} \geq 0$ of $u$ on $\mathbb{R}$. For $\mu-$a.e. $t\in \Omega, \Vert f \Vert_{(M_{1}+M_{2})(t,u)} = 0$ if and only if $f(t) = 0$ for an arbitrary $\lambda$. For $\mu-$a.e. $t\in \Omega$, for any scalar $\alpha, 
\Vert \alpha f \Vert_{(M_{1}+M_{2})(t,u)} = \vert \alpha \vert \Vert f \Vert_{(M_{1}+M_{2})(t,u)}$ of $u$ on $\mathbb{R}$; and since for $\mu-$a.e. $t\in \Omega,(M_{1}+M_{2})(t,u)$ is a convex function of $u$ on $\mathbb{R}$, let $ Sn_{(M_{1}+M_{2})} = \Vert f(t) \Vert_{(M_{1}+M_{2})} + \Vert g(t)\Vert_{(M_{1}+M_{2})}$, then
\begin{eqnarray*}
\int_{\Omega}(M_{1}+M_{2})\left( t, \frac{\Vert f(t)+g(t)\Vert_{BS}}{Sn_{(M_{1}+M_{2})}} \right) d\mu 
& \leq & \frac{\Vert f(t) \Vert_{(M_{1}+M_{2})}}{Sn_{(M_{1}+M_{2})}}\int_{\Omega} (M_{1}+M_{2})\left( t, \frac{\Vert f(t)\Vert_{BS}}{Sn_{(M_{1}+M_{2})}} \right) d\mu \\
&+& \frac{\Vert g(t) \Vert_{(M_{1}+M_{2})}}{Sn_{(M_{1}+M_{2})}}\int_{\Omega} (M_{1}+M_{2}) \left( t, \frac{\Vert g(t)\Vert_{BS}}{Sn_{(M_{1}+M_{2})}} \right) d\mu \\
& \leq & \frac{\Vert f(t) \Vert_{(M_{1}+M_{2})}}{Sn_{(M_{1}+M_{2})}}\int_{\Omega} (M_{1}+M_{2}) \left( t, \frac{\Vert f(t)\Vert_{BS}}{\Vert f(t) \Vert_{(M_{1}+M_{2})}} \right) d\mu \\
&+& \frac{\Vert g(t) \Vert_{(M_{1}+M_{2})}}{Sn_{(M_{1}+M_{2})}}\int_{\Omega} (M_{1}+M_{2}) \left( t, \frac{\Vert g(t)\Vert_{BS}}{\Vert g(t) \Vert_{(M_{1}+M_{2})}} \right) d\mu \\
& \leq & 1,
\end{eqnarray*}
because $\int_{\Omega} (M_{1}+M_{2}) \left( t, \frac{\Vert f(t)\Vert_{BS}}{\Vert f(t) \Vert_{(M_{1}+M_{2})}} \right) d\mu \leq 1$ and $\int_{\Omega} (M_{1}+M_{2}) \left( t, \frac{\Vert g(t)\Vert_{BS}}{\Vert g(t) \Vert_{(M_{1}+M_{2})}} \right) d\mu \leq 1$, so for $\mu-$a.e. $t\in \Omega, \Vert f(t)+g(t)\Vert_{(M_{1}+M_{2})(t,u)} \leq \Vert f(t)\Vert_{(M_{1}+M_{2})(t,u)} + \Vert g(t)\Vert_{(M_{1}+M_{2})(t,u)}$ of $u$ on $\mathbb{R}$; that is, for $\mu-$a.e. $t\in \Omega,\Vert f \Vert_{(M_{1}+M_{2})(t,u)}$ is norm on $L_{(M_{1}+M_{2})(t,u)}(\Omega,\Sigma,\mu)$ of $u$ on $\mathbb{R}$. Thus, for $r \geq 1$ that $\left( L_{rM}(\Omega,\Sigma,\mu), \Vert \cdot \Vert_{rM} \right)$ and $\left( L_{M_{1}+M_{2}}(\Omega,\Sigma,\mu), \Vert \cdot \Vert_{M_{1}+M_{2}} \right)$ are Musielak-Orlicz spaces generated by $rM$ and $M_{1}+M_{2}$ with the luxemburg norm respectively.
  
\begin{flushright}
$\square$
\end{flushright}
 
\paragraph{Remark 3.4} If $M_{i}:\Omega \times \mathbb{R} \to \mathbb{R}, i = 1,2$ are Musielak $N$-functions, then the subtraction in $M_{1}-M_{2}$ does not preserve the positivity and the convexity of the Musielak $N$-functions, so $M_{1}-M_{2}$ is not necessary to be a Musielak $N$-function and if so, it would be as $M_{1}+M_{2}$ and the $\left( L_{M_{1}-M_{1}}(\Omega,\Sigma,\mu), \Vert \cdot \Vert_{M_{1}-M_{1}} \right)$ would be as $\left( L_{M_{1}+M_{1}}(\Omega,\Sigma,\mu), \Vert \cdot \Vert_{M_{1}+M_{1}} \right)$. Moreover, $M_{1}+k, k \in \mathbb{R} \backslash \left\lbrace 0 \right\rbrace, M_{1}M_{2}, M_{1}^{n}, n \in \mathbb{N} $ and $M_{1}/M_{2}$ are not Musielak $N$-functions by theorem 2.2, where for $\mu-$a.e. $t\in \Omega$ that $(M_{1}+k)(t,0) \neq 0, (M_{1}M_{2})(t,0) \neq 0, M_{1}^{n}(t,0) \neq 0$ and $\frac{M_{1}(t,0)}{M_{2}(t,0)} \neq 0$ .

\paragraph{Theorem 3.10.} If $M:\Omega \times \mathbb{R} \to \mathbb{R}$ is a bounded Musielak $N$-function, then there exists a sequence of Musielak $N$-function $\{\varphi_{n}:n \in \mathbb{N}\}, \varphi_{n}:\Omega \times \mathbb{R} \to \mathbb{R}$ such that $\varphi_{n} \to M$ on $\Omega \times \mathbb{R}$.

\paragraph{Proof.} It is given that $M$ is a bounded Musielak $N$-function, so for each $u \in \mathbb{R}, M(t,u)$ is bounded and $\mu-$measurable function of $t$ on $\Omega$; so, there exists a sequence of simple functions $\{\varphi_{n}:n \in \mathbb{N}\}, \varphi_{n}:\Omega \times \mathbb{R} \to \mathbb{R}$ such that for each $u \in \mathbb{R}$, for all $\varepsilon > 0, \exists N \in \mathbb{N}, \vert M(t,u) - \varphi_{n}(t,u) \vert < \varepsilon$ for all $n \geq N$, for all $t\in \Omega$ by the basic approximation, then such convergence is uniform on $\Omega$ and pointwise on $\mathbb{R}$. Then, for all $n \in \mathbb{N}, \exists N \in \mathbb{N}$ such that for $\mu-$a.e. $t\in \Omega, \varphi_{n}(t,u)$ can satisfy the conditions $(1-4)$ of definition 2.1, and for each $u \in \mathbb{R}, \varphi_{n}(t,u)$ is $\mu-$measurable function of $t$ on $\Omega$ for all $n \geq N$, that is, these simple functions $\varphi_{n}, n \geq N$ are Musielak $N$-functions converge to $M$ on $\Omega \times \mathbb{R}$ as $n \to \infty$.

\begin{flushright}
$\square$
\end{flushright}

\paragraph{Theorem 3.11.} If $\left( L_{M}(\Omega,\Sigma,\mu), \Vert \cdot \Vert_{M} \right)$ is a Musielak-Orlicz space generated by a bounded Musielak $N$-function $M:\Omega \times \mathbb{R} \to \mathbb{R}$, then there exists a sequence of Musielak-Orlicz spaces $\left\lbrace \left( L_{\varphi_{n}}(\Omega,\Sigma,\mu), \Vert \cdot \Vert_{\varphi_{n}} \right), n \in \mathbb{N} \right\rbrace$, generated by a sequence of Musielak $N$-function $\{\varphi_{n}:n \in \mathbb{N}\}, \varphi_{n}:\Omega \times \mathbb{R} \to \mathbb{R}$ respectively, such that $\left( L_{\varphi_{n}}(\Omega,\Sigma,\mu), \Vert \cdot \Vert_{\varphi_{n}} \right) \to \left( L_{M}(\Omega,\Sigma,\mu), \Vert \cdot \Vert_{M} \right)$ via $\{\varphi_{n}:n \in \mathbb{N}\}$ as $n \to \infty$.

\paragraph{Proof.} It is given that $M:\Omega \times \mathbb{R} \to \mathbb{R}$ is a bounded Musielak $N$-function. By theorem 3.10, there exists a sequence of Musielak $N$-functions $\{\varphi_{n}:n \in \mathbb{N}\},\varphi_{n}:\Omega \times \mathbb{R} \to \mathbb{R}$ such that $\varphi_{n}$ converge to $M$ uniformly on $\Omega$ and pointwisely on $\mathbb{R}$ as $n \to \infty$; so $\{\varphi_{n}:n \in \mathbb{N}\}$ is uniformly bounded on $\Omega$. By the Lebesgue's dominated convergence theorem, for $\mu-$a.e. $t\in \Omega$,
\begin{eqnarray*}
L_{M(t,u)}(\Omega,\Sigma,\mu) &=& \left\lbrace f \in BS_{\Omega}: \int_{\Omega}M(t,\Vert f(t) \Vert_{BS})d\mu < \infty \right\rbrace \\
&=& \left\lbrace f \in BS_{\Omega}: \int_{\Omega} \lim \limits_{n \to \infty} \varphi_{n} \left( t, \Vert f(t) \Vert_{BS} \right) d\mu < \infty \right\rbrace \\
&=& \left\lbrace f \in BS_{\Omega}: \lim \limits_{n \to \infty} \int_{\Omega} \varphi_{n} \left( t, \Vert f(t) \Vert_{BS} \right) d\mu < \infty \right\rbrace \\
&=& \lim \limits_{n \to \infty} L_{\varphi_{n}(t,u)}(\Omega,\Sigma,\mu),
\end{eqnarray*}
and for $f \in L_{M(t,u)}(\Omega,\Sigma,\mu)$ that
\begin{eqnarray*}
\parallel f \parallel_{M(t,u)} &=& \inf \left\lbrace  \lambda > 0:\int_{\Omega} M \left(t, \frac{\Vert f(t) \Vert_{BS}}{\lambda}\right)d\mu \leq 1 \right\rbrace \\
&=& \inf \left\lbrace \lambda > 0: \int_{\Omega} \lim \limits_{n \to \infty} \varphi_{n} \left(t, \frac{\Vert f(t) \Vert_{BS}}{\lambda}\right) d\mu < 1 \right\rbrace \\
&=& \inf \left\lbrace \lambda > 0:\lim \limits_{n \to \infty} \int_{\Omega} \varphi_{n} \left(t, \frac{\Vert f(t) \Vert_{BS}}{\lambda}\right) d\mu < 1 \right\rbrace \\
&=& \lim \limits_{n \to \infty} \parallel f \parallel_{\varphi_{n}(t,u)};
\end{eqnarray*}
that is, for all $n \in \mathbb{N}$, for $\mu-$a.e. $t\in \Omega, \Vert f \Vert_{\varphi_{n}(t,u)}$ can satisfy the norm's properties on $L_{\varphi_{n}(t,u)}(\Omega,\Sigma,\mu)$ of $u$ on $\mathbb{R}$. Therefore, for all $n \in \mathbb{N}, \exists N \in \mathbb{N}, \left( L_{\varphi_{n}}(\Omega,\Sigma,\mu),\Vert \cdot \Vert_{\varphi_{n}}\right)$ is a Musielak-Orlicz spaces generated by a simple function $\varphi_{n}$ with the luxemburg norm for all $n \geq N$ respectively and $\left( L_{\varphi_{n}}(\Omega,\Sigma,\mu),\Vert \cdot \Vert_{\varphi_{n}}\right) \to \left( L_{M}(\Omega,\Sigma,\mu), \Vert \cdot \Vert_{M} \right)$ as $n \to \infty$.

\begin{flushright}
$\square$
\end{flushright}

\section{Properties, calculus and basic approximation of Musielak-Orlicz functions}
In this section we going to study properties, calculus and basic approximation of Musielak-Orlicz functions and the Musielak-Orlicz space generated by them and the relation between the Musielak $N$-function and the Musielak-Orlicz function and their Musielak-Orlicz spaces. Some theorems will be left without proof because they are a generalization to the ones that we did in section 3, so we can follow the similar way to proof them by consider the conditions of the Musielak-Orlicz function.

\paragraph{Theorem 4.1.} Let $M$ be a Musielak $N$-function and $MO$ be a Musielak-Orlicz function such that for $\mu-$a.e. $t\in \Omega$, 
\begin{align*}
\lim \limits_{u \to 0} \frac{MO(t,u)}{u} \neq 0, 
\end{align*}
then $L_{M}(\Omega,\Sigma,\mu) \neq L_{MO}(\Omega,\Sigma,\mu).$

\paragraph{Proof.} Assume $L_{M}(\Omega,\Sigma,\mu) = L_{MO}(\Omega,\Sigma,\mu)$, this means $L_{M}(\Omega,\Sigma,\mu) \subseteq L_{MO}(\Omega,\Sigma,\mu)$ and $L_{MO}(\Omega,\Sigma,\mu) \subseteq L_{M}(\Omega,\Sigma,\mu)$, so there exist two positive constants $r_{1}$ and $r_{2}$ and $u_{0} > 0$ such that for $f \in L_{M}(\Omega,\Sigma,\mu) \left( = L_{MO}(\Omega,\Sigma,\mu) \right)$, for $\mu-$a.e. $t\in \Omega$, 
\begin{align*}
r_{1} \int_{\Omega}MO(t,\Vert f(t) \Vert_{BS})d\mu \leq \int_{\Omega}M(t,\Vert f(t) \Vert_{BS})d\mu < \infty, \ f \neq 0
\end{align*}
and 
\begin{align*}
 r_{2} \int_{\Omega}M(t,\Vert f(t) \Vert_{BS})d\mu \leq \int_{\Omega}MO(t,\Vert f(t) \Vert_{BS})d\mu < \infty, \ f \neq 0.
\end{align*}
So for $\mu-$a.e. $t\in \Omega $,
\begin{align*}
r_{1} MO(t, u) \leq M(t, u), \ u \geq u_{0}
\end{align*}
and 
\begin{align*} 
 r_{2} M(t, u) \leq MO(t, u), \ u \geq u_{0}.
\end{align*}
Taking the limit $u \to 0$ we get for $\mu-$a.e. $t\in \Omega $, 
\begin{align*}
r_{1} \lim \limits_{u \to 0} \frac{MO(t, u)}{u} \leq \lim \limits_{u \to 0} \frac{M(t, u)}{u} = 0 
\end{align*}
and 
\begin{align*}
0 = r_{2} \lim \limits_{u \to 0} \frac{M(t, u)}{u} \leq \lim \limits_{u \to 0} \frac{MO(t, u)}{u}. 
\end{align*}
By the Squeez theorem for functions, we get for $\mu-$a.e. $t\in \Omega $,
\begin{align*}
\lim \limits_{u \to 0} \frac{MO(t, u)}{u} = 0 
\end{align*}
which contradicts the assumption.

\begin{flushright}
$\square$
\end{flushright}
 
\paragraph{Theorem 4.2.} Let $M$ be a Musielak $N$-function and $MO$ be a Musielak-Orlicz function such that for $\mu-$a.e. $t\in \Omega$, 
\begin{align*}
\lim \limits_{u \to \infty} \frac{MO(t,u)}{u} \neq \infty, 
\end{align*}
then $L_{M}(\Omega,\Sigma,\mu) \neq L_{MO}(\Omega,\Sigma,\mu).$

\paragraph{Proof.}  Following the similar way of theorem 4.1's proof and letting $u$ go to $\infty$, we get for $\mu-$a.e. $t\in \Omega $, 
\begin{align*}
r_{1} \lim \limits_{u \to \infty} \frac{MO(t, u)}{u} \leq \lim \limits_{u \to \infty} \frac{M(t, u)}{u} = \infty 
\end{align*}
and 
\begin{align*}
\infty = r_{2} \lim \limits_{u \to \infty} \frac{M(t, u)}{u} \leq \lim \limits_{u \to \infty} \frac{MO(t, u)}{u}. 
\end{align*}
By the Squeez theorem for functions, we get for $\mu-$a.e. $t\in \Omega $,
\begin{align*}
\lim \limits_{u \to \infty} \frac{MO(t, u)}{u} = \infty
\end{align*}
which contradicts the assumption.

\begin{flushright}
$\square$
\end{flushright}

\paragraph{Theorem 4.3.} If $MO_{1}:\Omega \times [0,\infty) \to [0,\infty)$ is a Musielak-Orlicz function such that $\lim \limits_{u \to 0} \frac{MO_{1}(t,u)}{u} \neq 0$ or $\lim \limits_{u \to \infty} \frac{MO_{1}(t,u)}{u} \neq \infty$ and $MO_{2}:\Omega \times [0,\infty) \to [0,\infty)$ is a function such that for $\mu-$a.e. $t\in \Omega, MO_{1}(t,u) = MO_{2}(t,u)$ of $u$ on $[0,\infty)$, then $MO_{2}$ is a Musielak-Orlicz function not Musielak $N-$function.

\paragraph{Proof.}  Since $MO_{1}$ is a Musielak-Orlicz function and for $\mu-$a.e. $t\in \Omega, MO_{1}(t,u) = MO_{2}(t,u)$ of $u$ on $[0,\infty)$, then for $\mu-$a.e. $t\in \Omega,MO_{2}(t,u)$ is convex function of $u$ on $[0,\infty)$; $MO_{2}(t,0) = MO_{1}(t,0) = 0$ and $MO_{2}(t,u) = MO_{1}(t,u) > 0$ for $u \neq 0$. Moreover, for each $u \in [0,\infty), MO_{2}(t,u)$ is $\mu-$ measurable function of $t$ on $\Omega$. So, $MO_{2}$ is a Musilak-Orlicz function on $\Omega \times [0,\infty)$. Now, assume that $MO_{2}$ is a Musielak $N-$function, then for $\mu-$a.e. $t\in \Omega$, 
\begin{align*}
0 \neq \lim \limits_{u \to 0} \frac{MO_{1}(t,u)}{u} = \lim \limits_{u \to 0} \frac{MO_{2}(t,u)}{u} = 0,
\end{align*}
or
\begin{align*}
\infty \neq \lim \limits_{u \to \infty} \frac{MO_{1}(t,u)}{u} = \lim \limits_{u \to \infty} \frac{MO_{2}(t,u)}{u} = \infty,
\end{align*}
which is a contradiction.

\begin{flushright}
$\square$
\end{flushright}

\paragraph{Theorem 4.4.} If $\left( L_{MO_{1}}(\Omega,\Sigma,\mu),\Vert \cdot \Vert_{M_{1}}\right)$ is a Musielak-Orlicz space generated by a Musielak-Orlicz function $M_{1}$ with the luxemburg norm and $MO_{2}:\Omega \times [0,\infty) \to [0,\infty)$ is a function such that for $\mu-$a.e. $t\in \Omega, M_{1}(t,u) = M_{2}(t,u)$ of $u$ on $\mathbb{R}$, then for $\mu-$a.e. $t\in \Omega, \left( L_{MO_{1}}(\Omega,\Sigma,\mu),\Vert \cdot \Vert_{MO_{1}}\right) = \left( L_{MO_{2}}(\Omega,\Sigma,\mu),\Vert \cdot \Vert_{MO_{2}}\right)$ of $u$ on $\mathbb{R}$ and hence $\left( L_{MO_{2}}(\Omega,\Sigma,\mu),\Vert \cdot \Vert_{MO_{2}}\right)$ is a Musielak-Orlicz space generated by $MO_{2}$ with the luxemburg norm.

\paragraph{Theorem 4.5.} If $\left\lbrace MO_{n}: n \in \mathbb{N} \right\rbrace$ is a sequence of Musielak-Orlicz functions $MO_{n}:\Omega \times [0,\infty) \to [0,\infty)$ such that for all $ n \in \mathbb{N}, \lim \limits_{u \to 0} \frac{MO_{n}(t,u)}{u} \neq 0$ or $\lim \limits_{u \to \infty} \frac{MO_{n}(t,u)}{u} \neq \infty$, then
\begin{align*}
\sup_{n \in \mathbb{N}} MO_{n}, \inf_{n \in \mathbb{N}} MO_{n}, \lim \limits_{n \to \infty} \sup  MO_{n}, \lim \limits_{n \to \infty} \inf MO_{n}
\end{align*}
are Musielak-Orlicz functions not Musielak $N-$ functions on $\Omega \times [0,\infty)$.

\paragraph{Proof.} Since for all $n \in \mathbb{N},MO_{n}$ is a Musielak-Orlicz function, then for $\mu-$a.e. $t\in \Omega, \sup_{n \in \mathbb{N}} MO_{n}(t,u)$ and $\inf_{n \in \mathbb{N}} MO_{n}(t,u)$ are convex of $u$ on $[0,\infty)$; for $\mu-$a.e. $t\in \Omega, \sup_{n \in \mathbb{N}} MO_{n}(t,0) = 0$ and $\inf_{n \in \mathbb{N}} MO_{n}(t,0) = 0$; and for $\mu-$a.e. $t\in \Omega, \sup_{n \in \mathbb{N}} MO_{n}(t,u) > 0$ and $\inf_{n \in \mathbb{N}} MO_{n}(t,u) > 0$ for $u \neq 0$. Moreover, for each $u \in [0,\infty),\sup_{n \in \mathbb{N}} MO_{n}(t,u)$ and $\inf_{n \in \mathbb{N}} MO_{n}(t,u)$ are $\mu-$measurable function of $t$ on $\Omega$. So, $\sup_{n \in \mathbb{N}} MO_{n}$ and $\inf_{n \in \mathbb{N}} MO_{n}$ are Musielak-Orlicz functions on $\Omega \times [0,\infty)$. Now, assume that $ \sup_{n \in \mathbb{N}} MO_{n}$ and $\inf_{n \in \mathbb{N}} MO_{n}$ are Musielak $N$-functions, then for $\mu-$a.e. $t\in \Omega$,
\begin{align*}
\lim \limits_{u \to 0} \frac{\sup_{n \in \mathbb{N}} MO_{n}(t,u)}{u} = 0 \Leftrightarrow \mathrm{for \ all} \ n \in \mathbb{N}, \lim \limits_{u \to 0} \frac{MO_{n}(t,u)}{u} = 0,  \\
\lim \limits_{u \to 0} \frac{\inf_{n \in \mathbb{N}} MO_{n}(t,u)}{u} = 0 \Leftrightarrow \mathrm{for \ all} \ n \in \mathbb{N}, \lim \limits_{u \to 0} \frac{MO_{n}(t,u)}{u} = 0, 
\end{align*}
or
\begin{align*}
\lim \limits_{u \to \infty} \frac{\sup_{n \in \mathbb{N}} MO_{n}(t,u)}{u} = \infty \Leftrightarrow \mathrm{for \ all} \ n \in \mathbb{N}, \lim \limits_{u \to \infty} \frac{MO_{n}(t,u)}{u} = \infty,  \\
\lim \limits_{u \to \infty} \frac{\inf_{n \in \mathbb{N}} MO_{n}(t,u)}{u} = \infty \Leftrightarrow \mathrm{for \ all} \ n \in \mathbb{N}, \lim \limits_{u \to \infty} \frac{MO_{n}(t,u)}{u} = \infty, 
\end{align*}
which is a contradiction to the assumptions. Also, for $\mu-$a.e. $t\in \Omega$, 
\begin{align*}
\lim \limits_{n \to \infty} \sup_{n \in \mathbb{N}}  MO_{n}(t,u) = \inf_{n \in \mathbb{N}} \sup_{k \geq n}MO_{k}(t,u) \\
\lim \limits_{n \to \infty} \inf_{n \in \mathbb{N}}  M_{n}(t,u) = \sup_{n \in \mathbb{N}} \inf_{k \geq n}M_{k}(t,u),
\end{align*}
of $u$ on $[0,\infty)$, it follows that $\lim \limits_{n \to \infty} \sup MO_{n}$ and $\lim \limits_{n \to \infty} \inf MO_{n}$ are Musielak-Orlicz functions not Musielak $N$-functions on $\Omega \times [0,\infty)$.

\begin{flushright}
$\square$
\end{flushright}

\paragraph{Theorem 4.6.}Let $\{MO_{n}: n \in \mathbb{N}\}$ be a sequence of Musielak-Orlicz functions $MO_{n}:\Omega \times [0,\infty) \to [0,\infty)$ such that for each $u \in [0,\infty), MO_{n}(t,u) \leq MO_{n+1}(t,u)$ on $\Omega$ for all $n \in \mathbb{N}$  . If $\{ \left( L_{MO_{n}}(\Omega,\Sigma,\mu),\Vert \cdot \Vert_{MO_{n}}\right): n \in \mathbb{N} \} $ is a sequence of Musielak-Orlicz spaces generated by $\{ MO_{n}: n \in \mathbb{N} \}$ with the luxemburg norm respectively, then $\sup_{n \in \mathbb{N}}( L_{M_{n}}(\Omega,\Sigma,\mu),\Vert \cdot \Vert_{M_{n}}) = ( L_{S}(\Omega,\Sigma,\mu),\Vert \cdot \Vert_{S})$, $\inf_{n \in \mathbb{N}}( L_{MO_{n}}(\Omega,\Sigma,\mu),\Vert \cdot \Vert_{MO_{n}}) = ( L_{I}(\Omega,\Sigma,\mu),\Vert \cdot \Vert_{I})$, $\lim \limits_{n \to \infty} \sup( L_{MO_{n}}(\Omega,\Sigma,\mu),\Vert \cdot \Vert_{MO_{n}}) = ( L_{LS}(\Omega,\Sigma,\mu),\Vert \cdot \Vert_{LS})$ and $\lim \limits_{n \to \infty} \inf(L_{MO_{n}}(\Omega,\Sigma,\mu),\Vert \cdot \Vert_{MO_{n}}) = ( L_{LI}(\Omega,\Sigma,\mu),\Vert \cdot \Vert_{LI})$ via $ \{ MO_{n}: n \in \mathbb{N} \}$, and hence $( L_{S}(\Omega,\Sigma,\mu),\Vert \cdot \Vert_{S})$, $( L_{I}(\Omega,\Sigma,\mu),\Vert \cdot \Vert_{I})$, $( L_{LS}(\Omega,\Sigma,\mu),\Vert \cdot \Vert_{LS})$ and $( L_{LI}(\Omega,\Sigma,\mu),\Vert \cdot \Vert_{LI})$ are Musielak-Orlicz spaces generated by 
\begin{align*}
S = \sup_{n \in \mathbb{N}} MO_{n}, I = \inf_{n \in \mathbb{N}} MO_{n}, LS = \lim \limits_{n \to \infty} \sup  MO_{n} \ \mathrm{and} \ LI = \lim \limits_{n \to \infty} \inf MO_{n}
\end{align*}
with the luxemburg norm respectively.

\paragraph{Theorem 4.7.} If $ \left\lbrace MO_{n}: n \in \mathbb{N}\right\rbrace $ is a seqence of Musielak-Orlicz functions $ MO_{n}:\Omega \times [0,\infty) \to [0,\infty)$ such that $\lim \limits_{u \to 0} \frac{MO_{n}(t,u)}{u} \neq 0$ or $\lim \limits_{u \to \infty} \frac{MO_{n}(t,u)}{u} \neq \infty$ and $MO_{n} \to MO, MO:\Omega \times [0,\infty) \to [0,\infty) $ pointwisely as $n \to \infty$, then $MO$ is a Musielak-Orlicz function not Musielak $N$-function and satisfy the $\bigtriangleup_{2}$-condition.

\paragraph{Proof.} Since the convergence of $MO_{n}$ to $MO$ is pointwisely, then for $\mu-$a.e. $t\in \Omega$ that
\begin{eqnarray*}
MO(t,u) = \lim \limits_{n \to \infty} \inf_{n \in \mathbb{N}}MO_{n}(t,u) = \lim \limits_{n \to \infty} \sup_{n \in \mathbb{N}}MO_{n}(t,u)
\end{eqnarray*}
of $u$ on $\mathbb{R}$. Now, assume that $MO$ is  a Musielak $N$-function, then for $\mu-$a.e. $t\in \Omega$, 
\begin{align*}
0 \neq \lim \limits_{n \to \infty} \lim \limits_{u \to 0} \frac{MO_{n}(t,u)}{u} = \lim \limits_{u \to 0} \lim \limits_{n \to \infty}\frac{MO_{n}(t,u)}{u} = \lim \limits_{u \to 0} \frac{MO(t,u)}{u} = 0,
\end{align*}
or 
\begin{align*}
\infty \neq \lim \limits_{n \to \infty} \lim \limits_{u \to 0} \frac{MO_{n}(t,u)}{u} = \lim \limits_{u \to 0} \lim \limits_{n \to \infty}\frac{MO_{n}(t,u)}{u} = \lim \limits_{u \to 0} \frac{MO(t,u)}{u} = \infty,
\end{align*}
which is a contradiction. So, $MO$ is a Musielak-Orlicz function not Musielak $N$-function on $\Omega \times [0,\infty)$ according to theorem 4.5; and the $\bigtriangleup_{2}$-condition is clear to satisfy.

\begin{flushright}
$\square$
\end{flushright}
 
\paragraph{Theorem 4.8.} Let $\{MO_{n}: n \in \mathbb{N}\}$ be a sequence of Musielak-Orlicz functions $MO_{n}:\Omega \times [0,\infty) \to [0,\infty)$ such that for each $u \in [0,\infty), MO_{n}(t,u) \to MO(t,u)$ on $\Omega$ as $n \to \infty$ and $\vert MO_{n}(t,u) \vert \leq G(t,u)$, where for each $u \in [0,\infty), G$ is absolutely integrable on $\Omega$. If $\left\lbrace \left( L_{M_{n}}(\Omega,\Sigma,\mu),\Vert \cdot \Vert_{M_{n}} \right): n \in \mathbb{N} \right\rbrace $ is a sequence of Musielak-Orlicz spaces generated by $ \left\lbrace MO_{n}: n \in \mathbb{N}\right\rbrace $ with the luxemburg norm, then $ ( L_{MO_{n}}(\Omega,\Sigma,\mu),\Vert \cdot \Vert_{MO_{n}}) \to (L_{MO}(\Omega,\Sigma,\mu),\Vert \cdot \Vert_{MO})$ via $ \{ MO_{n}: n \in \mathbb{N} \}$ as $n \to \infty$ and $\left( L_{MO}(\Omega,\Sigma,\mu),\Vert \cdot \Vert_{MO} \right)$ is a Musielak-Orlicz space generated by $MO$ with the luxemburg norm.

\paragraph{Corollary 4.1.} Under theorem 4.8's assumptions with $\vert MO_{n}(t,u) \vert \leq MO(t,u)$ and $MO_{n}$ satisfies the $\bigtriangleup_{2}$-condition for all $n \in \mathbb{N}$, if $f \in \left( L_{MO}(\Omega,\Sigma,\mu),\Vert \cdot \Vert_{MO} \right)$ then there exists a sequence of functions $\{f_{n}:n \in \mathbb{N}\}$  such that $f_{n} \in \left( L_{MO_{n}}(\Omega,\Sigma,\mu),\Vert \cdot \Vert_{MO_{n}} \right)$ for all $n \in \mathbb{N}$ and $f_{n} \underset{n \to \infty} \longrightarrow f$ under the luxemburg norm $\Vert \cdot \Vert_{MO}$.

\paragraph{Theorem 4.9.} If $MO:\Omega \times [0,\infty) \to [0,\infty)$ is a Musielak-Orlicz function such that $\lim \limits_{u \to 0} \frac{MO(t,u)}{u} \neq 0$ and $M:\Omega \times \mathbb{R} \to \mathbb{R}$ is a Musielak $N$-function, then for $r \in \mathbb{R_{+}}$ that $rMO$ and $MO+M$ are Musielak-Orlicz functions not Musielak $N$-functions.  

\paragraph{Proof.} It is given that $MO$ is a Musielak-Orlicz function and $M$ is a Musielak $N$-function, then both $MO$ and $M$ are non-negative convex functions on $[0,\infty)$. So, for $r \in \mathbb{R_{+}}$, for $\mu-$a.e. $t\in \Omega, rMO(t,u)$ and $(MO+M)(t,u)$ are convex of $u$ on $\mathbb{R}$; for $\mu-$a.e. $t\in \Omega, rMO(t,0) = 0$ and $(MO+M)(t,0) = 0$; and for $\mu-$a.e. $t\in \Omega, rMO(t,u) > 0$ and $(MO+M)(t,u) > 0$ for $u \neq 0$. Moreover, for each $u \in [0,\infty),rMO(t,u)$ and $(MO+M)(t,u)$ are $\mu-$measurable functions of $t$ on $\Omega$. So, $rMO$ and $MO+M$ are Musielak-Orlicz functions on $\Omega \times [0,\infty)$. Now, assume that $rMO$ and $MO+M$ are Musielak $N$-functions, then for $\mu-$a.e. $t\in \Omega$, 
\begin{align*}
0 = \lim \limits_{u \to 0} \frac{rMO(t,u)}{u} = r \lim \limits_{u \to 0} \frac{MO(t,u)}{u} \neq 0
\end{align*}
and
\begin{align*}
0 = \lim \limits_{u \to 0} \frac{(MO+M)(t,u)}{u} = \lim \limits_{u \to 0} \frac{MO(t,u)}{u}+\lim \limits_{u \to 0} \frac{M(t,u)}{u} \neq 0;
\end{align*}
so, both make a contradiction.

\begin{flushright}
$\square$
\end{flushright}

\paragraph{Theorem 4.10.}If $\left( L_{MO}(\Omega,\Sigma,\mu),\Vert \cdot \Vert_{MO} \right)$ is a Musielak-Orlicz space generated by a Musielak-Orlicz function $MO$ and $\left( L_{M}(\Omega,\Sigma,\mu),\Vert \cdot \Vert_{M} \right)$ generated by a Musielak $N$-function $M$ with the luxemburg norm, then $ \left( L_{rMO}(\Omega,\Sigma,\mu),\Vert \cdot \Vert_{rMO} \right), r \geq 1$ and $\left( L_{MO+M}(\Omega,\Sigma,\mu),\Vert \cdot \Vert_{MO+M} \right)$ are Musielak-Orlicz spaces generated by Musielak-Orlicz functions $rM$ and $MO+M$ with the luxemburg norm respectively.

\paragraph{Remark 4.1.} If $MO_{i} : \Omega \times [0,\infty) \to [0,\infty), i = 1,2$ are Musielak-Orlicz functions such that $\lim \limits_{u \to 0} \frac{MO_{i}(t,u)}{u} \neq 0,i = 1,2$ or $\lim \limits_{u \to \infty} \frac{MO_{i}(t,u)}{u} \neq \infty, i = 1,2$ and $M:\Omega \times \mathbb{R} \to \mathbb{R}$ is a Musielak $N$-function, then $MO_{1}+k, k \in \mathbb{R} \backslash \{0\}, MO_{1}MO_{2}, MO_{1}^{n}, n \in \mathbb{N}, \frac{MO_{1}}{MO_{2}},\frac{MO_{1}}{M}$ and $\frac{M}{MO_{1}}$ are neither Musielak-Orlicz functions nor Musielak $N$-functions because for $\mu-$a.e. $t\in \Omega, (MO_{1}+k)(t,0) \neq 0,(MO_{1}MO_{2})(t,0) \neq 0, MO_{1}^{n}(t,0) \neq 0, \frac{MO_{1}(t,0)}{MO_{2}(t,0)} \neq 0,\frac{MO_{1}(t,0)}{M(t,0)} \neq 0, \frac{M(t,0)}{MO_{1}(t,0)} \neq 0$. Also, $MO_{1} - MO_{2},MO_{1}-M$ and $M-MO_{1}$ are not necessary to be Musielak-Orlicz functions because the subtraction in them do not preserve the positivity and the convexity of the Musielak-Orlicz functions and if so, they would be as $MO_{1}+MO_{2}$ and $MO_{1}+M$ and $\left( L_{MO_{1}-MO_{2}}(\Omega,\Sigma,\mu),\Vert \cdot \Vert_{MO_{1}-MO_{2}} \right), \left( L_{MO_{1}-M}(\Omega,\Sigma,\mu),\Vert \cdot \Vert_{MO_{1}-M} \right) $ and $\left( L_{M-MO_{1}}(\Omega,\Sigma,\mu),\Vert \cdot \Vert_{M-MO_{1}} \right)$ would be as $\left( L_{MO_{1}+MO_{2}}(\Omega,\Sigma,\mu),\Vert \cdot \Vert_{MO_{1}+MO_{2}} \right)$ and $\left( L_{MO_{1}+M}(\Omega,\Sigma,\mu),\Vert \cdot \Vert_{MO_{1}+M} \right)$.

\paragraph{Theorem 4.11.} If $MO:\Omega \times [0,\infty) \to [0,\infty)$ is a bounded Musielak-Orlicz function, then there exists a sequence of Musielak-Orlicz function $\{\varphi_{n}:n \in \mathbb{N}\}, \varphi_{n}:\Omega \times [0,\infty) \to [0,\infty)$ such that $\lim \limits_{u \to 0} \frac{\varphi_{n}(t,u)}{u} \neq 0, \lim \limits_{u \to \infty} \frac{\varphi_{n}(t,u)}{u} \neq \infty$ and $\varphi_{n} \to M$ on $\Omega \times [0,\infty)$.

\paragraph{Proof.} It is given that $MO$ is a bounded Musielak-Orlicz function, so for each $u \in [0,\infty), MO(t,u)$ is bounded and $\mu-$measurable function of $t$ on $\Omega$; so there exists a sequence of simple functions $\{\varphi_{n}:n \in \mathbb{N}\}, \varphi_{n}:\Omega \times [0,\infty) \to [0,\infty)$ such that for each $u \in [0,\infty)$, for all $\varepsilon > 0, \exists N \in \mathbb{N}, \vert MO(t,u) - \varphi_{n}(t,u) \vert < \varepsilon$ for all $n \geq N$, for all $t\in \Omega$ by the basic approximation, then such convergence is uniform on $\Omega$ and pointwise on $[0,\infty)$. Then, for all $n \in \mathbb{N}, \exists N \in \mathbb{N}$ such that for $\mu-$a.e. $t\in \Omega, \varphi_{n}(t,u)$ can satisfy the conditions of Orlicz function on $[0,\infty)$, and for each $u \in [0,\infty), \varphi_{n}(t,u)$ is $\mu-$measurable function of $t$ on $\Omega$ for all $n \geq N$, then these simple functions $\varphi_{n}, n \geq N$ are Musielak-Orlicz functions converge to $MO$ on $\Omega \times [0,\infty)$ as $n \to \infty$. Now, assume $\{\varphi_{n}:n \in \mathbb{N}\}$ are Musielak $N$-functions, then

\begin{align*}
0 = \lim \limits_{n \to \infty}\lim \limits_{u \to 0} \frac{\varphi_{n}(t,u)}{u} =  \lim \limits_{u \to 0} \lim \limits_{n \to \infty} \frac{\varphi_{n}(t,u)}{u} = \lim \limits_{u \to 0} \frac{MO(t,u)}{u} \neq 0
\end{align*}
or
\begin{align*}
\infty = \lim \limits_{n \to \infty}\lim \limits_{u \to \infty} \frac{\varphi_{n}(t,u)}{u} =  \lim \limits_{u \to \infty} \lim \limits_{n \to \infty} \frac{\varphi_{n}(t,u)}{u} = \lim \limits_{u \to \infty} \frac{MO(t,u)}{u} \neq \infty
\end{align*}
which is a contradiction.

\begin{flushright}
$\square$
\end{flushright}

\paragraph{Theorem 4.12.} If $\left( L_{MO}(\Omega,\Sigma,\mu), \Vert \cdot \Vert_{MO} \right)$ is a bounded Musielak-Orlicz space generated by a Musielak-Orlicz function $MO:\Omega \times [0,\infty) \to [0,\infty)$, then there exists a sequence of Musielak-Orlicz spaces $\{( L_{\varphi_{n}}(\Omega,\Sigma,\mu), \Vert \cdot \Vert_{\varphi_{n}}), n \in \mathbb{N}\}$, generated by a sequence of Musielak-Orlicz functions $\varphi_{n}:\Omega \times [0,\infty) \to [0,\infty)$ respectively, such that $\left( L_{\varphi_{n}}(\Omega,\Sigma,\mu), \Vert \cdot \Vert_{\varphi_{n}} \right) \to \left( L_{MO}(\Omega,\Sigma,\mu), \Vert \cdot \Vert_{MO} \right)$ via $ \{ MO_{n}: n \in \mathbb{N} \}$ as $n \to \infty$.

\section{Examples}
\subsection{Examples of Musielak $N$-functions.}
\begin{enumerate}
\item  Every $N$-function is a Musielak $N$-function.
\item $M:\mathbb{R} \times \mathbb{R} \to [0,\infty),M(t,u) = (tu)^{2}$ is Musielak $N$-function, where
	  \begin{itemize}
	  \item for $\mu-$a.e. $t\in \mathbb{R}, M(t,u)$ is even convex because $M(t,-u) =  M(t,u)$ and $	\frac{\partial^{2} 		
	  		M(t,u)}{\partial u^{2}} = 2t^{2} \geq 0$ for all $u \in \mathbb{R}$
	  \item for $\mu-$a.e. $t\in \mathbb{R},M(t,u) = (tu)^{2} > 0$ for any $u > 0$
	  \item for $\mu-$a.e. $t\in \mathbb{R}, \lim \limits_{u \to 0} \frac{M(t,u)}{u} = \lim \limits_{u \to 0}{t^{2}u} = 0$
      \item for $\mu-$a.e. $t\in \mathbb{R}, \lim \limits_{u \to \infty} \frac{M(t,u)}{u} =  \lim \limits_{u \to \infty}{t^{2}	             u} = \infty$
	  \item for each $u \in \mathbb{R}, M(t,u) =(tu)^2$ is a $\mu-$measurable function of $t$ on $\mathbb{R}$  since it is 	 
	        continuous on measurable set $\mathbb{R}.$
	  \end{itemize}
\item $M:\mathbb{R} \times \mathbb{R} \to \mathbb{R}, M(t,u) = exp({\vert u \vert} + {\vert t \vert})-\vert u      
       \vert - \vert t \vert$ is Musielak $N$-function, where
	  \begin{itemize}
	  \item for $\mu-$a.e. $t\in \mathbb{R}, M(t,u)$ is even convex because $M(t,-u) = M(t,u)$ and $\frac{\partial^{2} M(t,u)}	
	        {\partial u^{2}} = exp({\vert u \vert} + {\vert t \vert}) > 0$ for all $u \in \mathbb{R}$ 
	   \item for $\mu-$a.e. $t\in \mathbb{R}, M(t,u) = exp({\vert u \vert} + {\vert t \vert})-\vert u \vert - \vert t \vert >  
	         0$ for any $u > 0$
	   \item for $\mu-$a.e. $t\in \mathbb{R}, \lim \limits_{u \to 0} \frac{M(t,u)}{u} = \lim \limits_{u \to 0} \frac{exp({\vert  
	           u 	
	   		 \vert}	+ {\vert t \vert})-\vert u \vert - \vert t \vert}{u} = 0$
	   \item for $\mu-$a.e. $t\in \mathbb{R}, \lim \limits_{u \to 0} \frac{M(t,u)}{u} = \lim \limits_{u \to \infty} 	 
	         \frac{exp({\vert u	\vert} + {\vert t \vert})-\vert u \vert - \vert t \vert}{u} = \infty$
	   \item for each $u \in \mathbb{R}, M(t,u) = exp({\vert u \vert} + {\vert t \vert})-\vert u \vert - \vert t \vert$ is a $	
	         \mu-$measurable function of $t$ on $\mathbb{R}$  since it is continuous on measurable set $\mathbb{R}.$
	  \end{itemize}	  
\end{enumerate}

\subsection{Examples of Musielak-Orlicz functions that are not Musielak $N$-functions.}
\begin{enumerate}
\item $M:[0,\infty)\times [0,\infty)\longrightarrow [0,\infty), M(t,u) = a^{tu}-1, a > 1$;
	\begin{itemize}
	\item for $\mu-$a.e. $t\in [0,\infty), a^{tu}-1 > 0$ on $(0,\infty)$ and $a^{tu}-1 = 0$ whenever $u=0$
	\item for $\mu-$a.e. $t\in [0,\infty), M(t,u) = a^{tu}-1$ is convex function on $\left[ 0,\infty \right)$ since 		 
		  $\frac{\partial^{2} M(t,u)}{\partial u^{2}}  = a^{tu}(t \log  a)^2 > 0$ for all $u \in [0,\infty)$
    \item for all $ u\geq 0, M \left(t,u\right) = a^{tu}-1$ is $\mu-$measurable function of $t$ on $[0,\infty)$ since it is 	 		 continuous on measurable set $[0,\infty)$
	\end{itemize} 
	then $M(t,u)$ is a Musielak-Orlicz function but not Musielak $N-$ function because for $\mu-$a.e. $t 	
	\in [0,\infty)$,
	$$ \lim \limits_{u \longrightarrow 0} \frac{a^{tu}-1}{u} = t \log a \neq 0 $$
\item $M:(- \infty,0)\times [0,\infty)\longrightarrow [0,\infty), M(t,u) = u^{t}$;
      \begin{itemize}
	   \item for $\mu-$a.e. $t\in (- \infty,0), u^{t} > 0$ on $(0,\infty)$ and $u^{t} = 0$ whenever
 		    $u=0$
 	   \item for $\mu-$a.e. $t\in (- \infty,0), M(t,u) = u^{t}$ is convex function on $\left[ 0,\infty \right)$ since 		 
		  $ \frac{\partial^{2} M(t,u)}{\partial u^{2}} = t(t-1)u^{t-2} > 0$ for all $u \in [0,\infty)$
		  \item for all $ u\geq 0, M \left(t,u\right) = u^{t}$ is $\mu-$measurable function of $t$ on $(- \infty,0)$ since it 	 				is continuous on measurable set $(- \infty,0)$
	    \end{itemize} 
	    then $M(t,u)$ is a Musielak-Orlicz function but not Musielak $N-$ function because for $\mu-$a.e. $t\in (- \infty,0)$,
	    $$ \lim \limits_{u \longrightarrow 0} \frac{u^{t}}{u} = \infty \neq 0 $$
\item $M:R\times [0,\infty)\longrightarrow [0,\infty), M(t,u) = (t+1)^{2}u$
	\begin{itemize}
	\item for $\mu-$a.e. $t\in R, M(t,u) > 0$ on $(0,\infty)$ and $M(t,u) = 0$ whenever 	
		  $u=0$
	\item for $\mu-$a.e. $t\in R, M(t,u) = (t+1)^{2}u$ is convex function on $(0,\infty)$ since for all 		  
		  $u_{1},u_{2} \in (0,\infty)$ and for all $ \lambda \in (0,1)$ that
			\begin{eqnarray*}
			 M(t,\lambda u_{1}+(1-\lambda)u_{2}) &=& (t+1)^{2}(\lambda u_{1}+(1-\lambda)u_{2}) \\
										 &\leq & \lambda(t+1)^{2} u_{1} + (1-\lambda)(t+1)^{2}u_{2} \\
										 & =& \lambda M(t,u_{1}) + (1-\lambda)M(t,u_{2}) \\
			\end{eqnarray*}	
	\item for all $ u\geq 0, M \left(t,u\right)= (t+1)^{2}u$ is $\mu-$measurable function of $t$ on $R$ since it is continuous 	
		  on measurable set $R$
	\end{itemize}
then $M(t,u)$ is a Musielak-Orlicz function but not Musielak $N-$ function because for $\mu-					$a.e. $t	\in R$,	
		$$ \lim \limits_{u \longrightarrow 0} \frac{(t+1)^{2}u}{u} = (t+1)^{2} \neq 0 $$
		$$ \lim \limits_{u \longrightarrow \infty} \frac{(t+1)^{2}u}{u} = (t+1)^{2} \neq \infty $$

\end{enumerate} 

\section{Conclusion}
The concept of $N$-function can be generalized to Musielak $N$-function as the concept of Orlicz function is generalized to Musielak-Orlicz function. $\mu-$almost everywhere property, supremum, infimum, limit, convergence and basic convergence of a sequence of Musielak $N-$functions and Musielak-Orlicz spaces generated by them can be considered using facts and results of the measure theory . Also, the relationship between Musielak $N$-functions and Musielak-Orlicz functions and Musielak-Orlicz spaces generated by them have been studied according to facts and results of the measure theory.

\end{document}